\documentclass{siamart220329}

\usepackage{amsmath,amssymb}
\usepackage{enumitem}
\usepackage{epsfig,graphics}
\usepackage{booktabs}
\usepackage{url}
\usepackage{tikz}
\usepackage{float}
\graphicspath{{./figures/}}

\usepackage{algorithm}
\usepackage{algpseudocode}

\usepackage{array,booktabs}

\usepackage{caption} 
\newdimen{\algindent}
\algnewcommand{\LeftComment}[1]{\Statex \hspace{\algindent}
\( \triangleright \) #1}

\renewcommand{\le}{\leqslant}

\def\beq{\begin{equation}}
\def\eeq{\end{equation}}
\def\eeq{\end{equation}}

\def\rprime{\hbox{\hskip.25em\raise.5ex\hbox{$'$}\hskip.15em}}

\def\ddn1{{\frac{\partial}{\partial \nu_{\yb}}}}

\newcommand{\ignore}[1]{}

\newtheorem{remark}[theorem]{Remark}

\newcommand{\be}{\begin{equation}}
\newcommand{\ee}{\end{equation}}
\newcommand{\ba}{\begin{aligned}}
\newcommand{\ea}{\end{aligned}}
\newcommand{\bea}{\begin{eqnarray}}
\newcommand{\eea}{\end{eqnarray}}

\def\ymsc{Yau Mathematical Sciences Center,
  Tsinghua University, Beijing China 100084}

\def\qiuzhen{Qiuzhen College, 
  Tsinghua University, Beijing China 100084}

\def\papertitle{A fast automatically adaptive integral equation solver for nonlinear two-point boundary value systems}

\title{\papertitle}
\author{
  Tianze Zhang\thanks{\qiuzhen\,
    ({\tt zhangtz21@mails.tsinghua.edu.cn}).} \and
 Yixuan Ma\thanks{Department of Mathematics,
   Tsinghua University, Beijing China 100084\,
    ({\tt mayx22@mails.tsinghua.edu.cn}).} \and
 Jun Wang\thanks{\ymsc\,
    ({\tt jwang2020@tsinghua.edu.cn}).} 
}

\makeatletter
\renewenvironment{proof}[1][Proof]{%
  \par\addvspace{1ex plus .25ex minus .1ex}%
  \noindent\textit{#1}.\ \ignorespaces
}{%
  \unskip\nobreak\hfill\proofbox\par\addvspace{1ex plus .25ex minus .1ex}%
}
\makeatother

\makeatletter
\def\cref@override@label@type#1{}
\makeatother

\begin{document}

\maketitle

\begin{abstract}
We present a fast, high-order, automatically adaptive integral-equation solver
for nonlinear two-point boundary value systems.  The method builds on the fast
direct solver for linear second-kind integral equations and completes the
mathematical formulation needed for general boundary conditions,
including cases in which the simplest background Green's function is singular
and cases with genuinely rank-deficient boundary conditions.  Nonlinear problems
are written as nonlinear second-kind integral equations and are solved by
combining this fast adaptive linear solver with Newton iteration and
user-specified continuation.  In this way the nonlinear solver inherits the full strengths of the linear
integral-equation method: a fast direct solver, arbitrary order accuracy, and
automatic adaptive resolution with only a small overhead over a
solve on the final mesh.  The resulting algorithm is particularly effective for problems with localized or
multiscale structure, including boundary layers, dense oscillations, transition
regions, hypersensitive optimal-control systems, and singular perturbations.
The same framework also provides a computational foundation for future studies
of bifurcation and solution-branch exploration.
\end{abstract}

\begin{keywords}
  nonlinear two-point boundary value problems, integral equations, fast direct solvers, adaptive refinement, continuation methods, singular perturbation
\end{keywords}

\begin{AMS}
31A10, 65F30, 65E05, 65Y20
\end{AMS}

\pagestyle{myheadings}
\thispagestyle{plain}
\markboth{T. Zhang, Y. Ma, J. Wang}
{Fast adaptive nonlinear BVP systems}

\tableofcontents 
\vspace{3mm}
\section{Introduction \label{sec:intro}}  

Many problems in science and engineering require the numerical solution of two-point boundary value problems for systems of ordinary differential equations, which take the form of
\begin{equation}  \label{bvpnonlin}
\begin{split}
\Phi'(x) &=  F(\Phi(x),x), \\
A\cdot \Phi(a) &+ C\cdot \Phi(c) =\gamma \;, 
\end{split}
\end{equation}
where $\Phi: [a,c]\rightarrow \mathbb{R}^n$ is in $C^1[a,c]$, $F:\mathbb{R}^n\times[a,c]\rightarrow \mathbb{R}^n$ is smooth, $A,C\in\mathbb{R}^{n\times n}$, and $\gamma\in\mathbb{R}^n$.  Such systems arise after rewriting high-order scalar equations, singular perturbation problems, optimal-control problems, and travelling-wave or bifurcation equations.  We focus on linear boundary conditions, which already cover a large class of applications.  More general nonlinear boundary conditions can in principle be handled by applying the same linearization idea to the boundary equations, but this extension is not pursued here.

The difficult regimes for two-point BVPs are often governed by local structure: boundary layers, internal transition regions, dense oscillations, and hypersensitive dynamics may occupy only a small part of the interval.  Direct discretization of the differential operator must then resolve these small scales while also solving the resulting linear systems inside Newton or continuation iterations.  Moreover, discretized differential operators inherit mesh-dependent conditioning; refining the mesh to resolve a layer can make the algebraic problem substantially more ill-conditioned than the underlying continuous BVP.  Continuation methods encounter the same issue repeatedly: as a parameter varies, layers and transition regions may appear or move, forcing standard differential-discretization approaches to combine path following with careful mesh adaptation and, in difficult regimes, problem-specific treatment \cite{newtconv,newt2,Allgower2003}.

Integral-equation formulations address this issue at its source.  By using a background Green's function, the BVP is converted to a second-kind integral equation, so the conditioning of the discrete system respects the conditioning of the underlying continuous problem rather than the mesh spacing.  Greengard and Rokhlin developed this idea for scalar two-point BVPs \cite{greengardrokhlin}, Lee and Greengard built an automatic adaptive solver for stiff scalar problems \cite{leegreengard}, and Leeb and Rokhlin later treated fourth-order scalar equations \cite{leebrokhlin2020}.  Starr and Rokhlin formulated and solved the first-order system case \cite{starrrokhlin}: Special structure of the Green's function gives a hierarchically off-diagonal low-rank representation, and the resulting continuous-level fast solver has linear complexity. This can be viewed as a special case of fast direct solvers,
where a compressed representation of the system's matrix is constructed once and then applied efficiently to one or many right-hand sides \cite{martinssonrokhlin2005}.

The system formulation of Starr and Rokhlin, however, imposed an artificial restriction in order to use the simplest background Green's function: the matrix $A+C$ had to be nonsingular.  This is not a solvability condition for the original BVP.  The authors proposed using a manually constructed linear transform to convert a user's boundary conditions into this form, but that requirement limits the black-box use of the algorithm.  More importantly, an arbitrary manual transform may introduce artificial ill conditioning and thereby compromise the principal advantage of the integral-equation formulation.

Our first contribution is to complete the mathematical formulation for general linear boundary conditions.  We give an explicit, well-conditioned-by-construction boundary transform whenever the boundary conditions are nondegenerate.  When the boundary conditions themselves are truly rank deficient, we first test compatibility and then complete the boundary conditions, reducing the problem to a family of full-rank problems.  In the linear case this yields the full affine solution space: a particular solution together with basis solutions for the missing boundary data.  Because the Starr--Rokhlin algorithm is a fast direct solver, this solution-space computation is naturally organized as a multiple-right-hand-side solve; the expensive factorization is shared and the extra basis functions are obtained at marginal cost.

Our second contribution is to extend this completed linear solver to nonlinear systems.  We formulate the nonlinear BVP as a nonlinear second-kind integral equation and use the fast adaptive linear solver as the computational kernel for Newton corrections and user-specified continuation.  In this way we have coupled the fast direct solver, boundary-condition preprocessing, Newton iteration, continuation, and automatic adaptivity into a nonlinear BVP solver which inherits the essential advantages of the linear integral-equation method: well-conditioned discretizations, arbitrary order accuracy, linear complexity, and automatic adaptivity.  These properties are especially valuable for nonlinear problems with localized or multiscale structure, where continuation and Newton iterations repeatedly require reliable linearized solves on evolving meshes.

The rest of the paper is organized as follows.  Section~\ref{sec:preli} develops the integral-equation formulation and gives the explicit boundary-condition preprocessing, including the treatment of truly degenerate boundary conditions.  Section~\ref{sec:fastlinear} summarizes the fast direct solver and high-order Nystr\"om discretization used for the linearized problems.  Sections~\ref{sec:newton} and~\ref{sec:continuation} describe Newton iteration and user-specified continuation for the nonlinear integral equation, and Section~\ref{sec:adaptivity} gives the adaptive strategy.  Section~\ref{sec:results} presents a sequence of challenging numerical examples designed to test the formulation, accuracy, scaling, adaptive efficiency, behavior on nonlinear problems with multiple solutions, and performance on singularly perturbed systems.  Section~\ref{sec:conclusion} concludes with the role of the method in future bifurcation computations.

\section{Integral formulation and boundary preprocessing}\label{sec:preli}

In this section, we develop the second-kind integral-equation formulation for two-point boundary value systems.  The linear formulation serves as the computational kernel for the nonlinear solver of later sections.
\subsection{Linear equations}
We start with definitions and notation.
\begin{definition}[Linear TPBVP]\label{linear-def}
    A linear TPBVP has the form
    \begin{gather}
        \Phi'(x) + p(x)\cdot \Phi(x) = f(x)\;, \label{eqn:lineqn} \\
        A\cdot \Phi(a) + C\cdot \Phi(c) = \gamma \;. \label{eqn:linbc}
    \end{gather}
    The equation is referred to as homogeneous if $f=0$:
    \begin{equation} \label{eqn:holineqn}
    \Phi'(x) + p(x)\cdot \Phi(x) = 0\;,
    \end{equation} 
    and the boundary condition is referred to as homogeneous if $\gamma=0$.
    \begin{equation}\label{eqn:holinbc}
    A\cdot \Phi(a) + C\cdot \Phi(c) =0 \;.
    \end{equation}
\end{definition}
\begin{definition}[Green's function]\label{greensfunc}
    A continuous function $G(x,t): [a,c]\times [a,c]\rightarrow \mathbf{L}(\mathbb{R}^{n\times n})$
    is referred to as the Green's function for the BVP problem \eqref{eqn:holinbc}\eqref{eqn:holineqn} if:
    \begin{itemize}
    \item $\frac{\partial G(x,t)}{\partial x}$ is continuous except at $x=t$.
    \item $G(x+0,x)-G(x-0,x)=I_n$ for all $x\in [a,c]$.
    \item $\frac{\partial}{\partial }G(x,t)+p(x)\cdot G(x,t)=0$ for all $x,t\in[a,c],\; x\neq t$.
    \item $A\cdot G(a,t)+C\cdot G(c,t)=0$ for all $t\in [a,c]$.
    \end{itemize}
\end{definition}
\begin{remark}
    Assume that $\Gamma(x)$ is the fundamental matrix of \eqref{eqn:holinbc}\eqref{eqn:holineqn}. If $D=A\cdot\Gamma(a)+C\cdot \Gamma(c) $ is nonsingular, then the Green's function exists and is unique~\cite{bliss1926}. Moreover, there is an unique solution $\Phi(x)$ for the inhomogeneous problem \eqref{eqn:holinbc} \eqref{eqn:lineqn}:
    \begin{equation}\label{sol_int}
        \Phi(x)=\int_a^c G(x,t) f(t)\;dt\;.
    \end{equation}   \end{remark}

For most variable-coefficient boundary value problems, the exact Green's function is unavailable.  We instead introduce a background problem for which the fundamental matrix and Green's function are known.  This converts the original equation with homogeneous boundary conditions into a second-kind integral equation.  The construction follows \cite{greengardrokhlin,starrrokhlin}.

\begin{definition}
    The background equation is a homogeneous equation with homogeneous boundary equations
\begin{equation} \label{background}
    \begin{split}
    \Phi'(x) + p_0(x)\cdot \Phi(x) &= 0\\
    A\cdot \Phi(a) + C\cdot \Phi(c) &=0 \;,
    \end{split}
\end{equation}
with the fundamental matrix $\Gamma_0(x)$ and Green's function $G_0$, satisfying $D_0=A\cdot\Gamma_0(a)+C\cdot\Gamma_0(c)$ is nonsingular.
\end{definition}
\begin{theorem}\label{thm:integ}
Consider the BVP
\begin{equation} \label{homobcbvp}
\begin{split}
\Phi'(x) + p(x)\cdot \Phi(x) &= f(x), \\
A\cdot \Phi(a) + C\cdot \Phi(c) &=0 .
\end{split}
\end{equation}
If it has a fundamental matrix $\Gamma(x)$ such that
$D:=A\cdot\Gamma(a)+C\cdot\Gamma(c)$ is nonsingular, then the unique solution can be represented as
\begin{equation}\label{intGreen}
    \Phi(x)=\int_a^c G_0(x,t)\sigma(t)dt,
\end{equation}
where $\sigma$ is the solution of the second-kind integral equation
\begin{equation}\label{integ}
\sigma(x)+(p(x)-p_0(x))\cdot \int_a^c G_0(x,t)\cdot \sigma(t)\;dt = f(x).
\end{equation}
\end{theorem}

The reduction of inhomogeneous boundary data to homogeneous form depends on the processed boundary condition matrices; it is deferred to Section~\ref{sec:inhomogeneous_bc_reduction} after the boundary-condition preprocessing.

\subsection{Background equations and the artificial restriction in the simplest formulation}
The simplest background equation is obtained by taking $p_0=0$. Then the Green's function is 
\begin{equation}\label{eqn:greensfunc1}
\overline{G}_0(x,t)=
\begin{cases}
I_n -(A+C)^{-1}\cdot C \;\;&(t\leq x)\;, \\
-(A+C)^{-1}\cdot C \;\; &(t\geq x) \;.
\end{cases}
\end{equation}
and Theorem~\ref{thm:integ} reduces to the following corollary.
\begin{corollary}\label{cor:greens0}
    If $D_0=A+C$ and $D=A\cdot\Gamma(a)+C\cdot\Gamma(c)$ is nonsingular, where $\Gamma(x)$ is the fundamental matrix of \eqref{eqn:holineqn} \eqref{eqn:holinbc}. Then the unique solution of \eqref{eqn:lineqn} \eqref{eqn:holinbc} is 
    \begin{equation}\label{intGreen0}
        \Phi(x)=\int_a^c \overline{G}_0(x,t)\sigma(t)dt,
    \end{equation}
    where $\sigma$ is the solution of the integral equation
    \begin{equation}\label{integ0}
    \sigma(x)+p(x)\cdot \int_a^c \overline{G}_0(x,t)\cdot \sigma(t)\;dt = f(x)\;,     
    \end{equation}
\end{corollary}

Two different issues must be separated.  The condition $\det(A+C)\neq0$ is only the condition that makes the particular background Green's function \eqref{eqn:greensfunc1} well defined; it is not a solvability condition for the original BVP.  Thus an otherwise well-posed problem may be excluded for a purely artificial reason.  This occurs, for example, when a second-order scalar equation with separated Dirichlet data is rewritten as a first-order system: the boundary matrices select the first component at opposite endpoints and $A+C$ is singular.  Higher-order scalar equations with separated boundary data lead to the same phenomenon.  In contrast, the condition $R(A)+R(C)\ne\mathbb R^n$ means that the boundary conditions themselves are rank deficient.  The former difficulty should be removed by a well-conditioned transform; the latter requires a compatibility check and a completion of the missing boundary data.

This distinction is essential for nonlinear problems.  Newton and continuation methods repeatedly solve linearized boundary value systems with the same boundary matrices.  A nonlinear solver built on the integral formulation is therefore only as general as its linear correction solver.  Completing the boundary-condition formulation is what makes the fast direct linear solver a reliable kernel for the nonlinear method.

\subsection{Automatic preprocessing for the case $\det(A+C)=0$}

We now give a constructive treatment of the case $\det(A+C)=0$.  This is not merely an implementation device; it completes the integral-equation formulation for general nondegenerate boundary conditions.  The background problem with $p_0=0$ is attractive because its Green's function is explicit, cheap to evaluate and leads to a well conditioned integral equation, but it is usable only when the transformed boundary matrix is invertible.  The following change of variables supplies this missing ingredient.

Let
\begin{equation}
    \Phi(x)=T(x)\hat\Phi(x),\qquad T(x)\in GL_n(\mathbb{R}),\quad x\in[a,c].
\end{equation}
The transformed equation and boundary condition are
\begin{align}
&\hat\Phi'(x)+T(x)^{-1}\cdot(T'(x)+p(x)\cdot T(x))\cdot\hat\Phi (x)=T(x)^{-1}f(x), \label{eqn:eqntrans}\\
&A\cdot T(a)\cdot\hat\Phi(a)+C\cdot T(c)\cdot\hat\Phi(c)=0. \label{eqn:bctrans}
\end{align}
Thus the transformed problem can be paired with the simplest background Green's function \eqref{eqn:greensfunc1} as soon as
\begin{equation}\label{eqn:trans_nonzero}
    \det(A T(a)+C T(c))\neq 0.
\end{equation}
The following theorem shows that, except for the genuinely rank-deficient boundary conditions discussed in the next subsection, such a transformation always exists and can be constructed explicitly.

\begin{theorem}[construction of an admissible boundary transform]\label{thm:bc_transform}
Let $A,C\in\mathbb{R}^{n\times n}$.  There exists a function $T\in C^1([a,c];\mathbb{R}^{n\times n})$ such that
\begin{equation}
    T(x)\in GL_n(\mathbb{R})\quad\text{for all }x\in[a,c],
    \qquad
    \det(A T(a)+C T(c))\neq 0,
\end{equation}
if and only if
\begin{equation}\label{eqn:gfun_cond2_draft}
    R(A)+R(C)=\mathbb{R}^n,
\end{equation}
or equivalently
\begin{equation}\label{eqn:gfun_cond1_draft}
    N(A^T)\cap N(C^T)=\{0\}.
\end{equation}
Moreover, $T$ may be chosen as a smooth path from the identity matrix to the product of a signed permutation matrix and a nonsingular diagonal matrix.
\end{theorem}

\begin{proof}
We first prove necessity.  If there exists a nonzero vector $z\in N(A^T)\cap N(C^T)$, then
\begin{equation}
    z^T(A T(a)+C T(c))=0
\end{equation}
for every choice of $T(a)$ and $T(c)$.  Hence $A T(a)+C T(c)$ cannot be nonsingular.  Therefore an admissible transform implies \eqref{eqn:gfun_cond1_draft}, which is equivalent to \eqref{eqn:gfun_cond2_draft} by elementary linear algebra.

For sufficiency, write the columns of $A$ and $C$ as
\begin{equation}
    A=(a_1,\ldots,a_n),\qquad C=(c_1,\ldots,c_n).
\end{equation}
Condition \eqref{eqn:gfun_cond2_draft} means that the $2n$ columns $a_1,\ldots,a_n,c_1,\ldots,c_n$ span $\mathbb{R}^n$.  Therefore there are index sets
\begin{equation}
    \{i_1,\ldots,i_r\},\qquad \{j_1,\ldots,j_{n-r}\}
\end{equation}
such that
\begin{equation}
    a_{i_1},\ldots,a_{i_r},c_{j_1},\ldots,c_{j_{n-r}}
\end{equation}
are linearly independent.  Choose a permutation $\sigma\in S_n$ such that
\begin{equation}
    \{i_1,\ldots,i_r\}\cup \{\sigma(j_1),\ldots,\sigma(j_{n-r})\}=\{1,\ldots,n\}.
\end{equation}
We look for a transform satisfying $T(a)=I$ and $T(c)=P_\sigma\Lambda$, where $P_\sigma$ is the permutation matrix associated with $\sigma$ and $\Lambda=\operatorname{diag}(\lambda_1,\ldots,\lambda_n)$ is diagonal.  Then
\begin{equation}
    A T(a)+C T(c)=A+C P_\sigma\Lambda.
\end{equation}
Choose a scalar $\lambda>1$ and set
\begin{equation}\label{eqn:lambdai_draft}
\left\{
\begin{aligned}
    &|\lambda_k|=\lambda, && k=\sigma(j_l),\quad l=1,\ldots,n-r,\\
    &|\lambda_k|=\lambda^{-1}, && k=i_l,\quad l=1,\ldots,r.
\end{aligned}
\right.
\end{equation}
The signs of the $\lambda_k$ may be chosen arbitrarily, or chosen to improve the conditioning of the resulting matrix.  With this convention,
\begin{equation}
    A+C P_\sigma\Lambda=(a_1+\beta_1 c_{\sigma(1)},\ldots,a_n+\beta_n c_{\sigma(n)}),
\end{equation}
where $\beta_k$ is the $k$th diagonal entry of $\Lambda$.  After scaling the columns indexed by $\sigma(J)$ by $\lambda^{-1}$ and the columns indexed by $I$ by $\lambda$, the limit as $\lambda\to\infty$ is a matrix whose columns are, up to signs, the independent set
\begin{equation}
    a_{i_1},\ldots,a_{i_r},c_{j_1},\ldots,c_{j_{n-r}}.
\end{equation}
Therefore $\det(A+C P_\sigma\Lambda)$ is not identically zero as a function of $\lambda$, and there exists a finite value of $\lambda$ for which \eqref{eqn:trans_nonzero} holds.

It remains to connect $I$ to $P_\sigma\Lambda$ inside $GL_n(\mathbb{R})$.  Decompose $P_\sigma$ into transpositions,
\begin{equation}
    P_\sigma=P_mP_{m-1}\cdots P_1.
\end{equation}
If $P_l$ swaps the coordinate directions $s_l$ and $k_l$, replace it by a smooth rotation path $\widetilde P_l(x)$ that is the identity on all other coordinates and has the two-dimensional block
\begin{equation}\label{eqn:rotation_path}
\begin{bmatrix}
\cos \theta(x) & -\sin \theta(x)\\
\sin \theta(x) & \phantom{-}\cos \theta(x)
\end{bmatrix},
\qquad
\theta(x)=\frac{\pi(x-a)}{2(c-a)}.
\end{equation}
At $x=c$, this realizes the transposition up to a sign change, which can be absorbed into the diagonal matrix.  Let
\begin{equation}
    \widetilde P_\sigma(x)=\widetilde P_m(x)\widetilde P_{m-1}(x)\cdots \widetilde P_1(x).
\end{equation}
For the diagonal part, choose a path
\begin{equation}\label{eqn:lambda_path}
    \Lambda(x)=\operatorname{diag}(\ell_1(x),\ldots,\ell_n(x)),
\end{equation}
where $\ell_k(a)=1$, $\ell_k(c)=\lambda_k$, and $\ell_k(x)$ never crosses zero.  For positive $\lambda_k$, one may simply take
\begin{equation}
    \ell_k(x)=1+\frac{x-a}{c-a}(\lambda_k-1).
\end{equation}
For negative $\lambda_k$, the sign can first be moved into the signed permutation part, leaving a positive diagonal path.  Finally set
\begin{equation}\label{eqn:T_path}
    T(x)=\widetilde P_\sigma(x)\Lambda(x).
\end{equation}
Then $T(a)=I$, $T(c)=P_\sigma\Lambda$ up to the prescribed signs, $T(x)$ is nonsingular for every $x$, and \eqref{eqn:trans_nonzero} holds.
\end{proof}

\begin{remark}
The theorem separates two issues that are often conflated.  The condition $\det(A+C)=0$ means only that the particular background Green's function \eqref{eqn:greensfunc1} is not directly available.  It does not indicate a defect in the boundary value problem.  In contrast, failure of \eqref{eqn:gfun_cond2_draft} means that the boundary condition itself is rank deficient for arbitrary boundary data.  That genuinely degenerate case is treated separately in Section~\ref{sec:true_degenerate_bc}.\end{remark}

\begin{remark}
The construction of $\Lambda$ in \eqref{eqn:lambdai_draft} uses a limiting argument.  Numerically, it is unnecessary to take $\lambda$ large in advance.  We start from $\lambda=1$ and double it until $A+C P_\sigma\Lambda$ is nonsingular and has acceptable conditioning.  In practice boundary conditions arising from physical problems are far from arbitrary, and the procedure above terminates with moderate values of $\lambda$ in all cases we have examined; a fully general bounded construction of $\Lambda$ is straightforward should it ever be needed.\end{remark}

The proof is constructive and gives the preprocessing algorithm below.  This step is performed once for a given pair $(A,C)$ and then reused throughout Newton iteration and continuation.

\begin{algorithm}[H]
\caption{Boundary-condition preprocessing for $\det(A+C)=0$}\label{alg:bc_preprocess}
\begin{algorithmic}[1]
\Require Boundary matrices $A=(a_1,\ldots,a_n)$ and $C=(c_1,\ldots,c_n)$ satisfying $R(A)+R(C)=\mathbb{R}^n$; interval endpoints $a,c$; conditioning threshold $\kappa_{\max}$.
\Ensure A $C^1$ path $T:[a,c]\to GL_n(\mathbb{R})$ such that $\det(A T(a)+C T(c))\neq0$.
\State Select index sets $I=\{i_1,\ldots,i_r\}$ and $J=\{j_1,\ldots,j_{n-r}\}$ such that $\{a_i:i\in I\}\cup\{c_j:j\in J\}$ is a basis of $\mathbb{R}^n$.
\State Choose $\sigma\in S_n$ such that $I\cup\sigma(J)=\{1,\ldots,n\}$.
\State Decompose $P_\sigma$ into transpositions $P_mP_{m-1}\cdots P_1$.
\State Replace each transposition $P_l$ by the rotation path \eqref{eqn:rotation_path}, absorbing endpoint signs into the diagonal matrix to obtain $\widetilde P_\sigma(x)$. 
\State Set $\lambda\gets1$.
\Repeat
    \State Form $\Lambda=\operatorname{diag}(\lambda_1,\ldots,\lambda_n)$ according to \eqref{eqn:lambdai_draft}.
    \State Set $D\gets A+C P_\sigma\Lambda$.
    \State If $D$ is singular or $\kappa_2(D)>\kappa_{\max}$, set $\lambda\gets2\lambda$.
\Until{$D$ is nonsingular and $\kappa_2(D)\leq\kappa_{\max}$}
\State Construct $\Lambda(x)$ by \eqref{eqn:lambda_path} and set $T(x)=\widetilde P_\sigma(x)\Lambda(x)$.
\State \Return $T$.
\end{algorithmic}
\end{algorithm}

After Algorithm~\ref{alg:bc_preprocess}, the transformed problem \eqref{eqn:eqntrans}--\eqref{eqn:bctrans} satisfies the nondegeneracy condition for the simplest background Green's function. The integral equation formulation \eqref{intGreen0}--\eqref{integ0} then readily applies.

\begin{remark}
The transform not only restores the black-box property, it also preserves conditioning.  Writing the linearized operator as $\mathcal L=\mathcal L_0+\mathcal E$, the second-kind equation contains $I+\mathcal E\mathcal L_0^{-1}$, so a poor transform or background Green's function can by itself degrade conditioning.  Since $T$, $T^{-1}$, and $T'$ are uniformly bounded, the original and transformed operators differ only by bounded pre- and post-composition; their condition numbers therefore differ by a constant depending on these bounds, which the construction keeps moderate.
\end{remark}

\subsection{Truly degenerate boundary conditions}\label{sec:true_degenerate_bc}

The discussion above assumes $R(A)+R(C)=\mathbb{R}^n$.  This is the regime in which arbitrary boundary data are consistent and the boundary-condition defect is only the artificial singularity of $A+C$ in the background Green's function.  There is also the genuinely rank-deficient case in which
\begin{equation}\label{eqn:true_degenerate_bc}
    R(A)+R(C)\neq \mathbb{R}^n.
\end{equation}
In that case the boundary condition itself does not impose $n$ independent conditions on arbitrary boundary data. The original problem may be inconsistent, or it may possess a family of solutions.  Examples of this category include incompletely specified endpoint data in optimal-control problems and certain free-boundary reductions.  We describe a boundary-completion procedure that reduces this case to the full-rank case treated above.

Let
\begin{equation}
    \mathcal B=[A\; C]\in\mathbb{R}^{n\times 2n},
    \qquad
    k=\operatorname{rank}(\mathcal B)=\dim(R(A)+R(C))<n.
\end{equation}
The boundary condition $A\Phi(a)+C\Phi(c)=\gamma$ is consistent only if
\begin{equation}\label{eqn:bc_compatibility}
    \gamma\in R(\mathcal B)=R(A)+R(C).
\end{equation}
Equivalently, if $Q_0\in\mathbb{R}^{n\times(n-k)}$ is an orthonormal basis of $R(\mathcal B)^\perp$, then
\begin{equation}\label{eqn:bc_compatibility_q0}
    Q_0^T\gamma=0.
\end{equation}
Let $Q_1\in\mathbb{R}^{n\times k}$ be an orthonormal basis of $R(\mathcal B)$.  Then the effective boundary conditions are precisely
\begin{equation}\label{eqn:effective_bc}
    Q_1^TA\Phi(a)+Q_1^TC\Phi(c)=Q_1^T\gamma.
\end{equation}

The missing $n-k$ boundary constraints may be supplied by an artificial completion.  Choose matrices $E_a,E_c\in\mathbb{R}^{(n-k)\times n}$ such that
\begin{equation}\label{eqn:completed_boundary_matrix}
    \bar{\mathcal B}
    :=
    \begin{bmatrix}
    Q_1^TA & Q_1^TC\\
    E_a & E_c
    \end{bmatrix}
    \in\mathbb{R}^{n\times 2n}
\end{equation}
has rank $n$.  Define
\begin{equation}\label{eqn:completed_A_C}
    \bar A=\begin{bmatrix}Q_1^TA\\ E_a\end{bmatrix},
    \qquad
    \bar C=\begin{bmatrix}Q_1^TC\\ E_c\end{bmatrix},
    \qquad
    \bar\gamma(\eta)=\begin{bmatrix}Q_1^T\gamma\\ \eta\end{bmatrix},
    \quad \eta\in\mathbb{R}^{n-k}.
\end{equation}
Then $R(\bar A)+R(\bar C)=\mathbb{R}^n$, and for each fixed artificial boundary datum $\eta$, the completed boundary condition
\begin{equation}\label{eqn:completed_bc}
    \bar A\Phi(a)+\bar C\Phi(c)=\bar\gamma(\eta)
\end{equation}
falls into the full-rank framework of the previous subsection.

\begin{theorem}[boundary completion for rank-deficient boundary conditions]\label{thm:true_degenerate_completion}
Assume $\operatorname{rank}[A\;C]=k<n$ and $\gamma\in R(A)+R(C)$.  Let $Q_1,Q_0$ and $(\bar A,\bar C,\bar\gamma(\eta))$ be defined as above, with $\operatorname{rank}\bar{\mathcal B}=n$.  Then a function $\Phi$ satisfies the original boundary condition
\begin{equation}
    A\Phi(a)+C\Phi(c)=\gamma
\end{equation}
if and only if there exists an artificial datum $\eta\in\mathbb{R}^{n-k}$ such that
\begin{equation}
    \bar A\Phi(a)+\bar C\Phi(c)=\bar\gamma(\eta).
\end{equation}
Consequently, if $\mathcal S$ denotes the solution set of the original boundary value problem and $\mathcal S_\eta$ denotes the solution set of the completed problem with artificial datum $\eta$, then
\begin{equation}\label{eqn:solution_family_completion}
    \mathcal S=\bigcup_{\eta\in\mathbb{R}^{n-k}}\mathcal S_\eta.
\end{equation}
In particular, choosing one value of $\eta$ selects one representative problem from the family, while varying $\eta$ parametrizes the boundary degrees of freedom left unspecified by the original rank-deficient boundary condition.
\end{theorem}

\begin{corollary}[multiple right-hand sides for the completed problem]\label{cor:multiple_rhs_degenerate_bc}
Consider the linear boundary value problem
\begin{equation}
    \Phi'(x)+p(x)\Phi(x)=f(x),
    \qquad
    A\Phi(a)+C\Phi(c)=\gamma,
\end{equation}
and assume that the compatibility condition \eqref{eqn:bc_compatibility} holds.  Let $r=n-k$ be the number of missing boundary conditions and assume that, for the completed boundary condition $(\bar A,\bar C)$, the completed BVP is uniquely solvable for each artificial datum $\eta\in\mathbb{R}^r$.  Let $\Phi_0$ solve
\begin{equation}
    \Phi_0'(x)+p(x)\Phi_0(x)=f(x),
    \qquad
    \bar A\Phi_0(a)+\bar C\Phi_0(c)=
    \begin{bmatrix}Q_1^T\gamma\\0\end{bmatrix},
\end{equation}
and, for $j=1,\ldots,r$, let $\Phi_j$ solve
\begin{equation}
    \Phi_j'(x)+p(x)\Phi_j(x)=0,
    \qquad
    \bar A\Phi_j(a)+\bar C\Phi_j(c)=
    \begin{bmatrix}0\\e_j\end{bmatrix},
\end{equation}
where $\{e_j\}_{j=1}^r$ is the standard basis of $\mathbb{R}^r$.  Then the solution family of the original rank-deficient boundary condition is
\begin{equation}\label{eqn:degenerate_bc_multiple_rhs_family}
    \Phi_\eta(x)=\Phi_0(x)+\sum_{j=1}^r \eta_j\Phi_j(x),
    \qquad \eta=(\eta_1,\ldots,\eta_r)^T\in\mathbb{R}^r.
\end{equation}
Equivalently, after one fast direct factorization for the completed boundary condition, the whole family is obtained by applying the solver to $r+1$ right-hand sides.
\end{corollary}

The multiple-right-hand-side form is not only a convenient mathematical
description; it is also the fast version of the algorithm.  The direct solver
constructs an operator-dependent compressed representation on the adaptive tree.
This representation depends on the coefficient matrix and on the completed
boundary condition, but not on the particular boundary-data right-hand side.
Thus the particular solve and the auxiliary homogeneous solves in
Corollary~\ref{cor:multiple_rhs_degenerate_bc} can share the same factorized
operator.  The additional right-hand sides are applied as a standard
multiple-right-hand-side direct solve, so the dominant factorization cost is paid
once and the marginal cost of the missing-boundary-data basis functions is small.  The implementation of this multiple-right-hand-side solve inside the fast algorithm is described in Section~\ref{sec:fastlinear}.

\begin{remark}
The artificial rows $(E_a,E_c)$ are not unique.  A stable default is obtained by an SVD row completion.  Let
\begin{equation}
    B_{\rm eff}=\begin{bmatrix}Q_1^TA&Q_1^TC\end{bmatrix}\in\mathbb{R}^{k\times 2n}
\end{equation}
and compute an orthonormal basis $W\in\mathbb{R}^{2n\times(2n-k)}$ for the orthogonal complement of the row space of $B_{\rm eff}$.  Choose $r=n-k$ rows from $W^T$ to form $E\in\mathbb{R}^{r\times 2n}$, preferably so that $\bar{\mathcal B}$ is well conditioned, and split $E=[E_a\;E_c]$.  This guarantees full boundary rank, although the completed BVP conditioning still depends on the differential equation through $\bar A\Gamma(a)+\bar C\Gamma(c)$.  The completed problem may still have $\det(\bar A+\bar C)=0$, but it satisfies $R(\bar A)+R(\bar C)=\mathbb{R}^n$ and can therefore be handled by Theorem~\ref{thm:bc_transform}.
\end{remark}

\begin{remark}[nonlinear rank-deficient boundary conditions]
For each fixed artificial
datum $\eta$, the completed boundary condition defines an ordinary full-rank
nonlinear BVP, and the Newton--fast-solver machinery described below applies
without change.  The present implementation solves a selected completed
representative, by default the case $\eta=0$, and reports the rank defect.  At
any fixed nonlinear iterate, however, the linearized correction problem is a
linear BVP, so the affine correction space associated with the missing boundary
data can be computed by the same multiple-right-hand-side machinery as in
Corollary~\ref{cor:multiple_rhs_degenerate_bc}.  Exploiting these affine
linearized spaces to trace or parameterize the full nonlinear solution manifold
is a natural subject for subsequent bifurcation and branch-following work, but
is not part of the solver interface studied in this paper.\end{remark}

The previous cases can be summarized as a single preprocessing stage for boundary conditions.

\begin{algorithm}[htbp]
\caption{Unified boundary-condition preprocessing}\label{alg:unified_bc_preprocess}
\begin{algorithmic}[1]
\Require Boundary matrices $A,C\in\mathbb{R}^{n\times n}$ and boundary data $\gamma\in\mathbb{R}^n$.
\Ensure A family of full-rank completed boundary conditions and, when needed, an admissible transform for the simplest background Green's function.
\State Form $\mathcal B=[A\;C]$ and compute $k=\operatorname{rank}(\mathcal B)$.
\If{$k<n$}
    \State Compute orthonormal bases $Q_1$ of $R(\mathcal B)$ and $Q_0$ of $R(\mathcal B)^\perp$.
    \State Check compatibility: if $Q_0^T\gamma\neq0$, report that the boundary condition is inconsistent.
    \State Set $B_{\rm eff}\gets[Q_1^TA\;Q_1^TC]$ and compute an SVD basis $W$ for $\operatorname{Row}(B_{\rm eff})^\perp\subset\mathbb{R}^{2n}$.
    \State Choose $r=n-k$ rows from $W^T$ to form $E=[E_a\;E_c]$, preferably minimizing the condition number of $\bar{\mathcal B}$.
    \State Set $(A,C,\gamma)\gets(\bar A,\bar C,\bar\gamma(\eta))$ with artificial datum $\eta\in\mathbb{R}^{n-k}$.
\EndIf
\If{$\det(A+C)\neq0$}
    \State Use the simplest background Green's function directly.
\Else
    \State Apply Algorithm~\ref{alg:bc_preprocess} to construct $T$ with $\det(A T(a)+C T(c))\neq0$.
    \State Transform the equation and boundary condition by \eqref{eqn:eqntrans}--\eqref{eqn:bctrans}.
\EndIf
\end{algorithmic}
\end{algorithm}

\subsection{Reduction of inhomogeneous boundary data}\label{sec:inhomogeneous_bc_reduction}

After the preprocessing steps above, including completion and the optional transformation $T$, the boundary condition has the form
\begin{equation}
    A\Phi(a)+C\Phi(c)=\gamma,
    \qquad \det(A+C)\neq0.
\end{equation}
The inhomogeneous boundary data can then be removed exactly as in the original Starr--Rokhlin formulation.  Set
\begin{equation}\label{eqn:boundary_constant_lift}
    \Phi_b=(A+C)^{-1}\gamma.
\end{equation}
Then $A\Phi_b+C\Phi_b=\gamma$, and writing $\Phi=\tilde\Phi+\Phi_b$ reduces the boundary condition to
\begin{equation}
    A\tilde\Phi(a)+C\tilde\Phi(c)=0.
\end{equation}
For a linear equation $\Phi'+p(x)\Phi=f(x)$, this only changes the right-hand side to
\begin{equation}
    \tilde f(x)=f(x)-p(x)\Phi_b.
\end{equation}
For a nonlinear equation $\Phi'=F(\Phi,x)$, one obtains the homogeneous-boundary problem
\begin{equation}
    \tilde\Phi'=F(\tilde\Phi+\Phi_b,x),
    \qquad
    A\tilde\Phi(a)+C\tilde\Phi(c)=0.
\end{equation}

\subsection{Nonlinear equations}
We now turn to the nonlinear TPBVP in \eqref{bvpnonlin}.  In view of Section~\ref{sec:inhomogeneous_bc_reduction}, it suffices to state the homogeneous-boundary case.  A similar integral equation can then be established.

\begin{corollary}\label{cor:nongreens0}
Assume that the fixed linear boundary conditions have been preprocessed as above
and that the inhomogeneous boundary data have been lifted out, so that the
simplest background Green's function $\overline G_0$ in \eqref{eqn:greensfunc1}
is available.  If $F$ is smooth, then $\Phi$ solves the nonlinear equation
\begin{equation}  \label{eqn:hononlin}
\begin{split}
\Phi'(x) &=  F(\Phi(x),x), \\
A\cdot \Phi(a) &+ C\cdot \Phi(c) =0 \;, 
\end{split}
\end{equation}
if and only if
\begin{equation}\label{eqn:intGreennonlin}
    \Phi(x)=\int_a^c \overline{G}_0(x,t)\sigma(t)dt\;,
\end{equation}
where $\sigma$ is the solution of the integral equation
\begin{equation}\label{eqn:intgenonlin}
\sigma(x)= F\left(\int_a^c \overline{G}_0(x,t)\sigma(t)dt,x\right)\;.    
\end{equation}
\end{corollary}

The important point is that the unknown density $\sigma$ satisfies a nonlinear second-kind equation: identity plus a compact nonlinear operator.  This is the formulation on which Newton iteration, continuation, and adaptivity are built.

\section{Fast linear solver and high-order discretization}\label{sec:fastlinear}

Every Newton correction and every continuation correction in this paper reduces to a linear second-kind integral equation of the form
\begin{equation}\label{eqn:linear_kernel}
    \xi(x)+q(x)\int_a^c G_0(x,t)\xi(t)\,dt = r(x),
\end{equation}
where $q$ is either the coefficient of a linear problem or the Jacobian of the nonlinear right-hand side along the current iterate.  We solve \eqref{eqn:linear_kernel} by the fast direct solver of Starr and Rokhlin for systems of two-point boundary value problems \cite{starrrokhlin}, together with the adaptive high-order discretization philosophy of Lee and Greengard \cite{leegreengard}.  We summarize the parts of the construction needed for the present nonlinear solver.

Let
\begin{equation}
    P[\xi](x)=\xi(x)+q(x)\int_a^cG_0(x,t)\xi(t)\,dt.
\end{equation}
Away from the diagonal the Green's function can be written in rank-$n$ factorized form $G_0(x,t)=\gamma_l(x)\gamma_r(t)^T$.  Hence, if the interval is split into two adjacent subintervals, the off-diagonal restrictions of $P$ have rank at most $n$, while the two diagonal restrictions are again integral operators of the same type on the child intervals.  Thus the matrix obtained after discretization has a hierarchically off-diagonal low-rank structure, which can be exploited directly at the operator level.

Consider first a single interval bisected into two halves.  A block elimination of the resulting integral equation, together with the Sherman--Morrison--Woodbury formula for its rank-$n$ off-diagonal blocks, shows that the problem decouples: the unknown on each half satisfies a local integral equation of exactly the same form as \eqref{eqn:linear_kernel}, and the two local solutions are then joined through a small set of coupling coefficients.  The global solve on the interval therefore reduces to two independent problems of the same type, at half the size.

Applying this reduction recursively yields the fast algorithm.  Let the interval $[a,c]$ be partitioned into $N$ subintervals, obtained for simplicity by recursive bisection, so that they are the leaves of a binary tree.  On each leaf the local integral equation is solved independently.  The coupling coefficients obey parent--child recursion relations: they are obtained by computing a few inner products of the local solutions at the leaves, followed by one upward and one downward sweep over the tree, after which the density is recovered on every leaf from its local solution and its coefficients.  Concrete formulas are given in \cite{starrrokhlin}.

The same recursions extend to multiple right-hand sides.  In the local solve stage the factorization of each leaf matrix is shared and simply applied to every right-hand side, so the dominant cost is unchanged; in the coefficient-passing stage the coupling vectors are augmented to matrices whose extra columns carry the right-hand sides and are advanced by the same parent--child relations.  This gives a fast direct solver for multiple right-hand sides at minimal extra cost---precisely what the completion of Section~\ref{sec:true_degenerate_bc} requires.

We now describe the discretization of the local problems in this hierarchy.  It is high order and naturally compatible with adaptivity.  On each leaf interval we use a $p$-point Chebyshev representation and a Nystr\"om discretization based on Clenshaw--Curtis quadrature.  Enforcing the local integral equation at the $p$ Chebyshev nodes and applying the Clenshaw--Curtis rule to the integral term yields the local dense system
\begin{equation}\label{eq:local_system}
    (\mathcal I_{np}+\mathcal K_I)\,\xi_I=r_I,
    \qquad
    \mathcal K_I\in\mathbb{R}^{np\times np},
\end{equation}
where $\mathcal K_I$ is the quadrature approximation of the Green's-function operator on $I$.  For fixed $n$ and $p$, both the data stored per tree node and the cost of combining two children are bounded independently of the number of active nodes; consequently the total work scales as $O(N n^3 p)$ and storage as $O(N n^2 p)$ for $N$ discretization nodes.  Increasing $p$ gives arbitrary order accuracy for smooth data, and the adaptive strategy of Section~\ref{sec:adaptivity} changes $N$ by refining only the intervals where the local error indicator is large.

Finally, the dominant cost is the solution of the local equations at the leaves, and these equations are decoupled.  After a refinement or coarsening step only the modified leaves require new local solves, while the rest of the tree is reused unchanged; adaptivity is therefore naturally inexpensive.

\section{Newton's iteration}\label{sec:newton}
Assume that $F:\mathbb{R}^{n}\times[a,c]\rightarrow\mathbb{R}^n$ is continuously differentiable with respect to its first argument.  After the boundary preprocessing and the reduction to homogeneous boundary conditions, the nonlinear BVP is written in density form as
\begin{equation}\label{eq:integ-nonli}
    \mathcal N(\sigma)(x)
    :=
    \sigma(x)-F((K\sigma)(x),x)=0,
    \qquad
    (K\sigma)(x)=\int_a^cG_0(x,t)\sigma(t)\,dt.
\end{equation}
The associated solution is $\Phi=K\sigma$.  The Fr\'echet derivative of $\mathcal N$ at $\sigma_i$ is
\begin{equation}\label{eq:newton_derivative_density}
    D\mathcal N(\sigma_i)\xi
    =
    \xi(x)-D_1F(\Phi_i(x),x)(K\xi)(x),
    \qquad
    \Phi_i=K\sigma_i.
\end{equation}
Thus Newton's method computes the density correction $\xi_i$ from
\begin{equation}\label{eq:integ-newt}
    \xi_i(x)-D_1F(\Phi_i(x),x)\int_a^cG_0(x,t)\xi_i(t)\,dt
    =F(\Phi_i(x),x)-\sigma_i(x),
\end{equation}
and updates
\begin{equation}\label{eq:newton_update_density}
    \sigma_{i+1}=\sigma_i+\xi_i,
    \qquad
    \Phi_{i+1}=K\sigma_{i+1}.
\end{equation}
Equation \eqref{eq:integ-newt} is again a linear second-kind integral equation of the form \eqref{eqn:linear_kernel}, with $q(x)=-D_1F(\Phi_i(x),x)$, where $D_1F$ denotes the Jacobian of $F$ with respect to its first argument.  It is solved by the fast direct solver of Section~\ref{sec:fastlinear}.

Equivalently, with $\delta\Phi_i=K\xi_i$, equation \eqref{eq:integ-newt} is the Newton correction for the differential residual,
\begin{equation}\label{eq:newt_ode_equivalent}
\delta\Phi_i'(x)-D_1F(\Phi_i(x),x)\delta\Phi_i(x)
    =-\Phi_i'(x)+F(\Phi_i(x),x),
\qquad
A\delta\Phi_i(a)+C\delta\Phi_i(c)=0.
\end{equation}
The numerical method nonetheless solves the density equation \eqref{eq:integ-newt}: on each leaf it is discretized by the same local system \eqref{eq:local_system}, now with $q=-D_1F(\Phi_i(\cdot),\cdot)$, rather than by discretizing the differential operator \eqref{eq:newt_ode_equivalent} directly.

The standard Newton--Kantorovich theory applies to \eqref{eq:integ-nonli}; see, for example, \cite{newtconv,newt2,Kantor}.  In particular, if $\widetilde\sigma$ is a solution, $D\mathcal N(\widetilde\sigma)$ is invertible, and the initial density $\sigma_0$ is sufficiently close to $\widetilde\sigma$, then the iterates remain in a neighborhood of $\widetilde\sigma$ and converge quadratically.  In practice, this local nature of Newton's method motivates the continuation strategy described in the next section.

\section{Continuation strategy}\label{sec:continuation}

Newton's method is a local method: its fast convergence is obtained once the iterate lies in the basin of attraction of the desired solution.  We therefore use continuation as a practical outer strategy: a difficult target problem is reached through a sequence of easier nearby BVPs, and the solution from one stage supplies the initial guess for the next.  Let
\begin{equation}
    \mathcal N(\sigma)=0
\end{equation}
denote the target nonlinear integral equation.  A user-specified homotopy is a family of nonlinear equations
\begin{equation}\label{eq:homotopy_general}
    \mathcal H(\sigma,s)=0,
    \qquad 0\le s\le 1,
\end{equation}
such that the problem at $s=0$ is easy to solve and
\begin{equation}
    \mathcal H(\sigma,1)=\mathcal N(\sigma).
\end{equation}
Typical choices include continuation in a physical parameter, continuation in boundary data, continuation in the forcing term, or continuation in the strength of the nonlinear term.

The local mathematical justification is the standard implicit function theorem; see, e.g., \cite{newtconv,Allgower2003}.  If there exists a smooth solution path $\sigma(s)$ and the linearized operator $D_\sigma\mathcal H(\sigma(s),s)$ is invertible along the path, then the solution can be continued locally in $s$.  At each step, a predictor provides an initial guess at the next parameter value and Newton's method is used as a corrector.  Each Newton correction has the same second-kind integral form as \eqref{eq:integ-newt} and is solved by the fast adaptive direct linear solver described in Section~\ref{sec:fastlinear}.

For boundary value problems, continuation often creates or reveals multiscale structure: as $s$ varies, boundary layers, transition layers, and oscillatory regimes may appear even when the starting problem is benign \cite{newtconv,Allgower2003}.  In collocation or finite-difference methods this couples the continuation step size strongly to mesh adaptation.  In the integral formulation every correction equation remains a second-kind equation and the adaptive tree refines and coarsens only where necessary, so the mesh tracks the evolving solution without accumulating refinement; the mechanism is described in Section~\ref{sec:adaptivity}.

In the implementation used for the numerical examples, the continuation schedule is supplied by the user.  If $0=s_0<s_1<\cdots<s_m=1$, then the solution at $s_j$ is used as the initial guess for Newton iteration at $s_{j+1}$.  A simple version of the procedure is given in Algorithm~\ref{alg:continuation_basic}.  

\begin{algorithm}[H]
\caption{User-specified continuation with Newton correction}\label{alg:continuation_basic}
\begin{algorithmic}[1]
\Require Homotopy $\mathcal H(\sigma,s)$, parameter values $0=s_0<s_1<\cdots<s_m=1$, initial solution $\sigma_0$ at $s_0$.
\For{$j=0,1,\ldots,m-1$}
    \State Set the initial guess at $s_{j+1}$ to be the solution at $s_j$.
    \Repeat
        \State Form the Newton correction equation for $\mathcal H(\sigma,s_{j+1})=0$.
        \State Solve the resulting linear second-kind integral equation by the fast adaptive solver.
        \State Update the density and the corresponding solution.
    \Until{the nonlinear residual or Newton correction is below tolerance}
\EndFor
\State \Return the solution at $s_m=1$.
\end{algorithmic}
\end{algorithm}

The fast adaptive integral-equation solver described here is the computational kernel that carries out the Newton and continuation corrections.  Its well-conditioned formulation, linear complexity, arbitrary order accuracy, and automatic adaptivity are properties of the integral-equation substrate, and they are inherited without change by any outer continuation or bifurcation algorithm built on top of it: pseudo-arclength continuation, fold detection, nullspace computation, and branch switching add outer control logic without altering the kernel or its cost.  The kernel and the outer continuation framework therefore operate at distinct levels, making the present contribution orthogonal to, and directly compatible with, techniques such as those in AUTO.

\section{Strategy for adaptivity}\label{sec:adaptivity}

The adaptive strategy follows the philosophy of Lee and Greengard \cite{leegreengard} but is extended to the nonlinear case: local polynomial representations are inspected a posteriori, and the interval tree is refined or coarsened until local solutions on all leaf intervals are resolved.  On each leaf interval $I$, each component of the density or solution is represented in a Chebyshev basis of order $p$,
\begin{equation}
    u_j(x)|_I\approx \sum_{m=0}^{p-1} a_{j,m}^{I}T_m(x),
    \qquad j=1,\ldots,n.
\end{equation}
For smooth functions, the decay of the Chebyshev coefficients provides a sharp and inexpensive resolution indicator; see, for example, the approximation theory viewpoint in \cite{Trefethen2013}.  In the implementation we use the Lee--Greengard tail quantity componentwise.  If
\begin{equation}
    u_j(x)|_I\approx \sum_{m=0}^{p-1} a_{j,m}^{I}T_m(x),
\end{equation}
then the monitor for component $j$ on interval $I$ is
\begin{equation}\label{eq:cheb_tail_indicator}
    S_I^{(j)}=|a_{j,p-2}^{I}|+|a_{j,p-1}^{I}-a_{j,p-3}^{I}|.
\end{equation}
For systems, refinement is triggered if any component exceeds its threshold, and coarsening is accepted when all component monitors are small.

The tail monitor is more than a refinement heuristic: the following theorem shows that resolvedness of the density implies accuracy.

\begin{theorem}[resolvedness implies accuracy, after \cite{leegreengard}]\label{thm:resolved}
Represent the density on each leaf $I$ by a Chebyshev expansion of order $p$ with coefficients $a_{m}^{I}$, and suppose that, for every leaf $I$,
\begin{equation}
    |a_{p-2}^{I}|<\tfrac{\varepsilon}{2}\|\tilde f\|,
    \qquad
    |a_{p-1}^{I}|<\tfrac{\varepsilon}{2}\|\tilde f\|,
    \qquad
    \|f-\tilde f\|<\varepsilon\|\tilde f\|.
\end{equation}
Then the computed density $\widehat\sigma$ satisfies the relative error bound
\begin{equation}\label{eq:lg_error_bound}
    \frac{\|\sigma-\widehat\sigma\|}{\|\sigma\|}
    \;\le\;
    \kappa(P)\,\bigl(\|P-\widehat P\|+1\bigr)\,\varepsilon,
\end{equation}
where $\kappa(P)$ is the condition number of the integral operator $P$.
\end{theorem}

\begin{remark}
Theorem~\ref{thm:resolved} shows that local resolvedness is a sufficient condition for accuracy, so the tail monitor of \eqref{eq:cheb_tail_indicator} inherits the same status.  The estimate uses the fact that $P$ is a compact perturbation of the identity; it is not readily available when the differential operator is discretized directly, because that operator is unbounded.  The same perturbation argument applies componentwise to systems and yields the analogous bound.
\end{remark}

For linear problems, this is the direct system analogue of the Lee--Greengard adaptivity criterion.  One could in principle maintain separate trees for different components of the vector solution, but this only changes the number of discretization points by a constant factor.  For simplicity, we use a single tree and require all components to be resolved on every leaf.  If one of the component monitors is too large, the interval is bisected.  Coarsening is applied by merging sibling leaves when the merged representation is also resolved by the same tail criterion.  The result may retain a small amount of over-refinement, but it avoids oscillatory refine--coarsen behavior and has been robust in the examples.

For nonlinear problems, adaptivity is inherited from the linear correction solver but is applied to the nonlinear state.  At each Newton step, and at each continuation parameter value, the terminal mesh of the previous solve is used as the initial mesh for the next correction equation.  Since the Newton update is $\sigma_{\rm new}=\sigma_{\rm old}+\delta\sigma$, the quantity that must be resolved is the updated density $\sigma_{\rm new}$, not the correction alone.  The monitors are therefore computed from $\sigma_{\rm new}$ using the same componentwise Lee--Greengard tail criterion, and refinement and coarsening are applied by the two-sided rule above.  The outer adaptivity control stops when two consecutive resolved approximations agree in the mixed relative norm
\begin{equation}
    \varrho=\frac{\|u_{\rm prev}-u\|_2}{\|u_{\rm prev}+u\|_2},
    \label{eq:adapt_mixed_norm}
\end{equation}
with a fallback to $\|u_{\rm prev}-u\|_2$ when the denominator is small.  In contrast with the final mesh-doubling check used in some purely linear adaptive solvers, the nonlinear correction solve does not automatically double the terminal mesh after local convergence, which would propagate unnecessary refinement through Newton and continuation steps.  Instead, the local Chebyshev tail monitors drive the mesh as the nonlinear iteration proceeds.

The resulting strategy is automatic: no a priori information is supplied by the user.  The mesh is concentrated in physically active regions and remains coarse where the solution is smooth or nearly constant.  A schematic version is given in Algorithm~\ref{alg:adaptive_refinement}.

\begin{algorithm}[H]
\caption{Adaptive refinement by Chebyshev tail estimates}\label{alg:adaptive_refinement}
\begin{algorithmic}[1]
\Require Current interval tree, Chebyshev order $p$, tolerance $\varepsilon$.
\Repeat
    \State Solve the current correction equation on the active tree.
    \State Form the quantity to be resolved: the density for a linear problem, or $\sigma_{\rm new}=\sigma_{\rm old}+\delta\sigma$ for a nonlinear correction.
    \ForAll{leaf intervals $I$}
        \State Compute Chebyshev coefficients and component monitors $S_I^{(j)}$ from \eqref{eq:cheb_tail_indicator}.
        \If{any component monitor exceeds its threshold}
            \State Mark $I$ for bisection.
        \EndIf
    \EndFor
    \State Bisect all marked leaves and update the tree data.
    \State Merge sibling leaves whose combined representation is resolved by the same criterion.
\Until{the mixed relative change \eqref{eq:adapt_mixed_norm} is below $\varepsilon$ and no leaf is marked for refinement}
\end{algorithmic}
\end{algorithm}

\section{Numerical results}\label{sec:results}
We now illustrate the performance of the method on a sequence of examples designed to test different aspects of the solver: adaptive resolution of boundary layers, robustness on long intervals, continuation through multiple solution branches, and singularly perturbed nonlinear behavior.  Unless otherwise stated, the discretization order is fixed and the adaptive tree is generated automatically.

\subsection{Shock-wave boundary layer and adaptive fast scaling}\label{sec:ex_shock_wave}

The first example is the viscous shock-wave problem used by Lee and Greengard
\cite{leegreengard} to illustrate adaptive refinement for a stiff scalar BVP,
\begin{equation}\label{eq:shock_scalar}
    \varepsilon u''(x)+2x u'(x)=0,\qquad -1\le x\le 1,
    \qquad u(-1)=-1,\quad u(1)=1,
\end{equation}
with $\varepsilon=10^{-5}$.  The exact solution is
\begin{equation}
    u(x)=\frac{\operatorname{erf}(x/\sqrt{\varepsilon})}
    {\operatorname{erf}(1/\sqrt{\varepsilon})}.
\end{equation}
Introducing the usual first-order variables
\begin{equation}
    \phi_1=u,
    \qquad
    \phi_2=u',
\end{equation}
we obtain the system
\begin{equation}\label{eq:shock_system}
    \phi_1'=\phi_2,
    \qquad
    \phi_2'=-\frac{2x}{\varepsilon}\phi_2,
    \qquad
    \phi_1(-1)=-1,
    \quad
    \phi_1(1)=1.
\end{equation}
This separated Dirichlet boundary condition is a simple but important case in
which the matrix $A+C$ for the simplest background Green's function is singular.
The preprocessing step constructs an admissible boundary transform automatically.
For this example the transform is especially benign: sampling $T(x)$ on the
interval gives $\max_x \kappa_2(T(x))=1$ in our diagnostic, so the preprocessing
removes the artificial singularity of $A+C$ without degrading the conditioning of
the problem.

We first test the fixed-mesh direct solver.  Starting from the two intervals
$[-1,0]$ and $[0,1]$, we manually bisect the two intervals adjacent to the layer
at $x=0$ and feed the resulting boundary-layer mesh to the nonadaptive system solver.  The
number of bisections ranges from $3$ to $7$.  Table~\ref{tab:shock_scaling}
reports the total number of discretization nodes $N$, the total time $T(N)$, the
fraction of time spent in the local leaf solves, and the maximum error against
the exact solution.

\begin{table}[H]
\centering

\begin{tabular}{rrrrrr}
\toprule
bisections & intervals & $N$ & $T(N)$ (s) & local fraction & max error \\
\midrule
3 & 8  & 128 & 0.013466 & 0.489 & $1.64\times 10^{2}$ \\
4 & 10 & 160 & 0.013345 & 0.598 & $2.43\times 10^{1}$ \\
5 & 12 & 192 & 0.015501 & 0.613 & $5.14\times 10^{-1}$ \\
6 & 14 & 224 & 0.016577 & 0.674 & $3.76\times 10^{-3}$ \\
7 & 16 & 256 & 0.018132 & 0.699 & $1.78\times 10^{-6}$ \\
\bottomrule
\end{tabular}
\caption{Fixed-mesh nonadaptive shock-wave solves.  The local fraction is included explicitly because this stage is the dominant cost as the mesh is refined.}
\label{tab:shock_scaling}
\end{table}

\begin{figure}[H]
\centering
\includegraphics[width=0.68\textwidth]{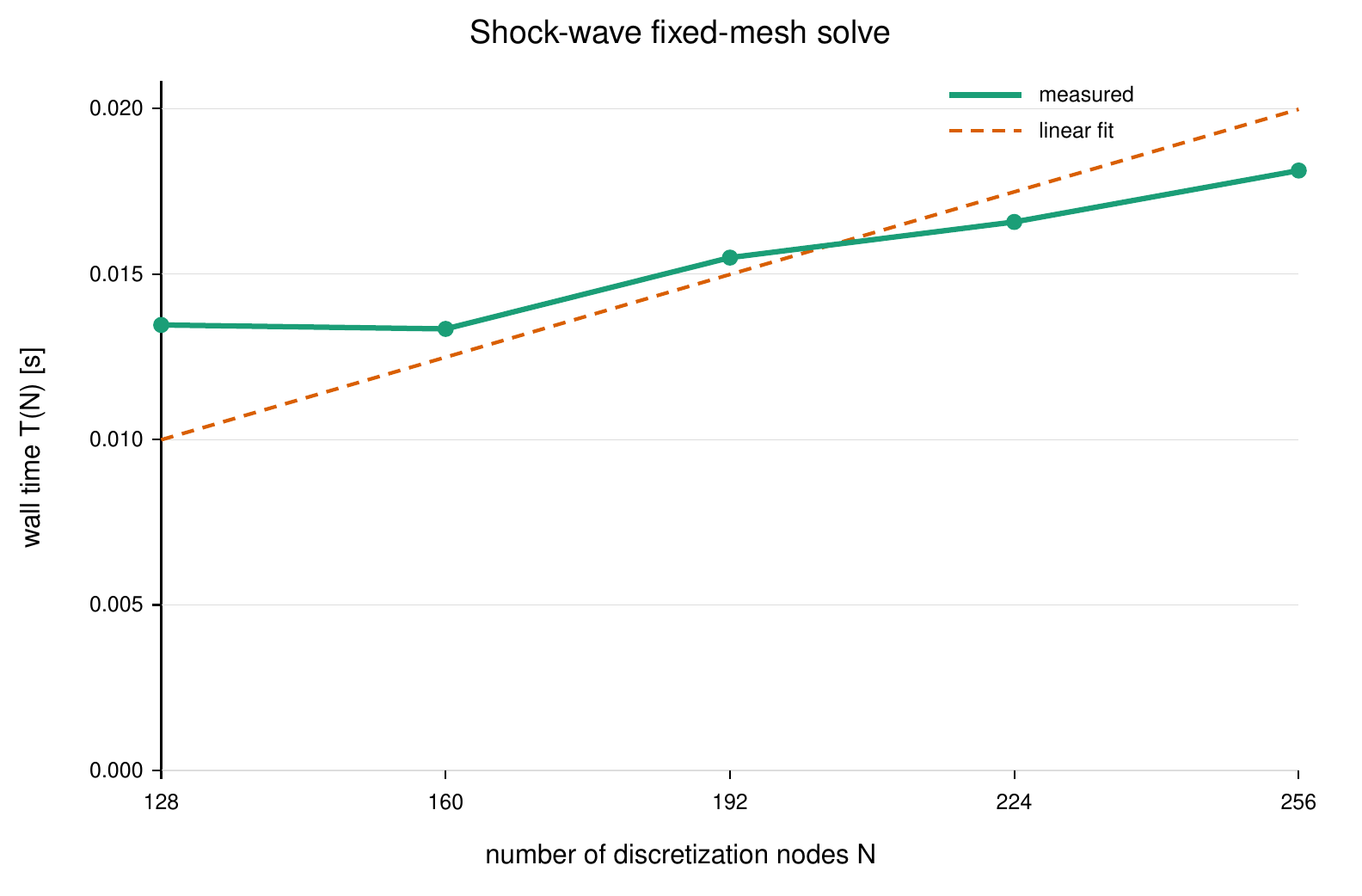}
\caption{Fixed-mesh shock-wave solve time versus the number of discretization
nodes.  The dashed curve is a least-squares linear fit through the origin.}
\label{fig:shock_scaling}
\end{figure}

Figure~\ref{fig:shock_scaling} plots $T(N)$ against $N$.  The measured time is
approximately linear over this range, as expected for the direct solver with
fixed system size and fixed Chebyshev order.

We next test the adaptive solver.  Starting from only the two initial intervals
$[-1,0]$ and $[0,1]$, the adaptive run produces a final mesh with $28$ intervals
and maximum error $1.32\times 10^{-10}$.  The adaptive run takes $6.25\times
10^{-2}$ seconds.  Running the nonadaptive solver once on this final adaptive
mesh takes $2.78\times10^{-2}$ seconds and gives the same error.  Thus the ratio
between the fully adaptive calculation and the direct solve on the final mesh is
approximately
\begin{equation}
    \frac{T_{\rm adaptive}}{T_{\rm final\ mesh}} \approx 2.29.
\end{equation}
This is the expected behavior for an efficient adaptive strategy: the code finds
the resolving mesh automatically at a cost only a small constant factor larger
than a single solve on the final mesh, while requiring no a priori knowledge of
the layer width.

\subsection{A truly degenerate boundary condition}\label{sec:ex_degenerate_bc}

The next example isolates the completion procedure for genuinely rank-deficient
boundary conditions while keeping the differential equation elementary.  Consider
the four-dimensional homogeneous system
\begin{equation}
\begin{aligned}
    u_1'&=u_2,        & u_2'&=-u_1,\\
    u_3'&=u_4,        & u_4'&=-u_3,
\end{aligned}
\qquad 0\le x\le 1,
\end{equation}
with only two independent boundary conditions,
\begin{equation}
    u_1(0)=0,
    \qquad
    u_3(0)=1.
\end{equation}
Thus $\operatorname{rank}[A\;C]=2<4$.  The exact solution set is the
two-parameter affine family
\begin{equation}
\begin{aligned}
    u_1(x)&=\eta_1\sin x,        & u_2(x)&=\eta_1\cos x,\\
    u_3(x)&=\cos x+\eta_2\sin x,& u_4(x)&=-\sin x+\eta_2\cos x,
\end{aligned}
\qquad \eta_1,\eta_2\in\mathbb R .
\end{equation}

The full-space interface first completes the boundary condition, solves the
particular problem and the two auxiliary homogeneous problems using a shared
multiple-right-hand-side direct solve, and returns a particular solution together
with two basis functions.  We compare this with a non-shared implementation
which applies the completed boundary conditions and calls the single-RHS solver
independently for the particular datum and the two auxiliary data.  The results
are shown in Table~\ref{tab:degenerate_bc}.  Both approaches recover the exact
affine family to roundoff accuracy, but the shared multiple-RHS version is almost
five times faster for this small test.

\begin{table}[H]
\centering
\begin{tabular}{lr}
\toprule
quantity & value \\
\midrule
nullity of boundary condition & 2 \\
fast full-space time & $8.88\times10^{-3}$ s \\
independent single-RHS time & $4.22\times10^{-2}$ s \\
time ratio & 4.76 \\
solution-space error & $4.44\times10^{-16}$ \\
free-component basis determinant & 1 \\
\bottomrule
\end{tabular}
\caption{Truly degenerate boundary condition example.  The fast full-space
interface shares the direct-solver factorization across the particular and
auxiliary right-hand sides.}
\label{tab:degenerate_bc}
\end{table}

This confirms the practical value of the completion formulation in
Section~\ref{sec:true_degenerate_bc}: the solver can automatically process a
rank-deficient boundary condition, return the affine solution space in the linear
case, and exploit the direct solver's multiple-right-hand-side structure rather
than solving each completed problem independently.

\subsection{Nonlinear mass-spring optimal control}\label{sec:ex_spring}

The third example is a nonlinear optimal-control problem for a spring-mass
system.  It tests two features of the solver which are absent from the preceding
linear examples: the use of Newton corrections for a genuinely nonlinear system,
and automatic mesh refinement for a solution with sharp endpoint transitions on
long time intervals.  The spring force is
\begin{equation}
    F_s(x)=-x-x^3,
\end{equation}
where $x=0$ is the unstretched position, and the scalar control force is denoted
by $u(t)$.  The control problem is to minimize
\begin{equation}
    J=\frac12\int_0^T (x^2+v^2+u^2)\,dt
\end{equation}
subject to
\begin{equation}
\left\{
    \begin{aligned}
        x'&=v,\\
        v'&=-x-x^3+u,\\
        x(0)&=1,\quad v(0)=0,\quad x(T)=0.75,\quad v(T)=0.
    \end{aligned}
\right.
\end{equation}
With adjoint variables $\lambda$ and $\mu$, the Hamiltonian is
\begin{equation}
    H=\frac12(x^2+v^2+u^2)+\lambda v+
      \mu(-x-x^3+u).
\end{equation}
Eliminating the control by the optimality condition gives the first-order BVP
\begin{equation}\label{eq:spring_optimality}
\left\{
    \begin{aligned}
        x'&=v,\\
        v'&=-x-x^3+u,\\
        \lambda'&=-x+\mu(1+3x^2),\\
        \mu'&=-v-\lambda,\\
        x(0)&=1,\quad v(0)=0,\quad x(T)=0.75,\quad v(T)=0.
    \end{aligned}
\right.
\end{equation}
This is the system solved by the code, with the control recovered from the
optimality condition inside the right-hand side evaluation.

\begin{figure}[t]
\includegraphics[width=1\textwidth]{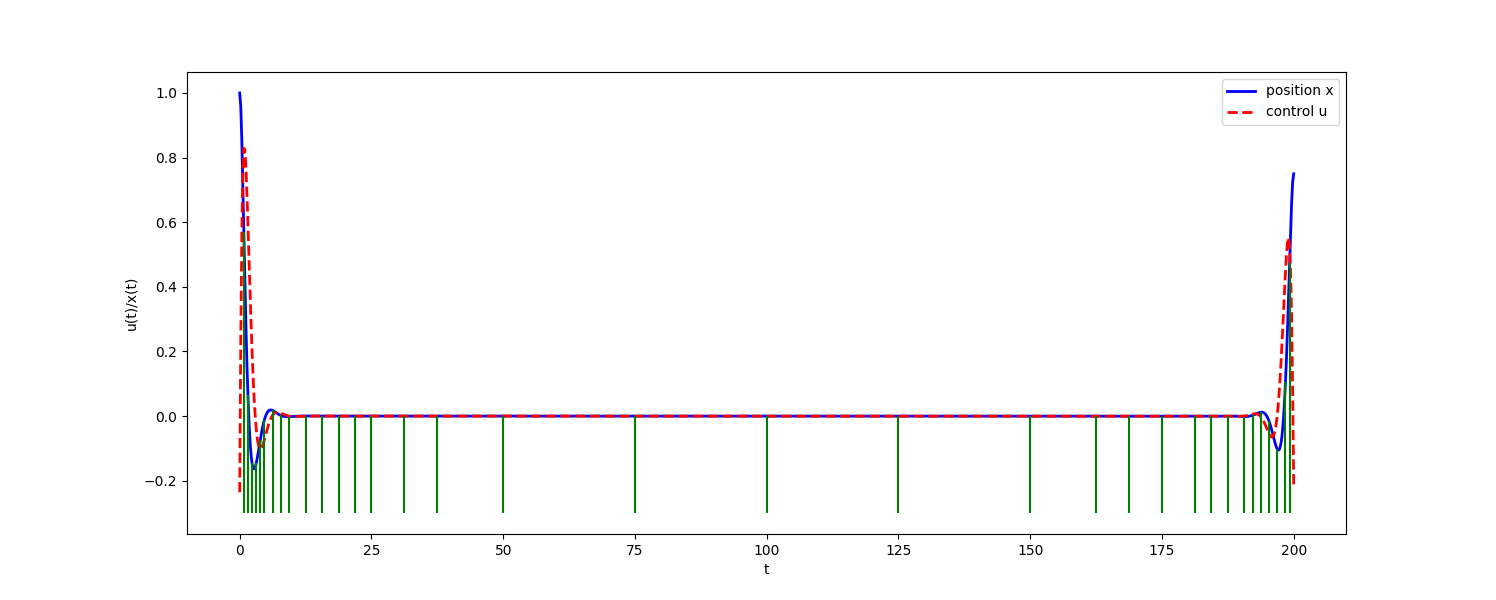}
\caption{Solution of the nonlinear optimal-control BVP
\eqref{eq:spring_optimality}. Green vertical lines indicate adaptive
subintervals.}
\label{fig:1-1}
\end{figure}

The solution develops rapid transitions near the endpoints while remaining
nearly flat over much of a long interval.  We use order $p=16$ and the simple
initial guess $x\equiv v\equiv\lambda\equiv\mu\equiv0$.  The resulting solution
and adaptive mesh are shown in Figure~\ref{fig:1-1}.  For $T=200$, the final mesh
contains 30 subintervals, with refinement concentrated near the endpoint
transition regions.

The same algorithm also applies on much longer time intervals.  We tested
intervals from $[0,20]$ to $[0,2.0\times10^5]$; mesh counts and run times are
reported in Table~\ref{tab:1}.  The mesh count grows slowly compared with the
length of the interval because most of the interior remains smooth and nearly
constant.
\begin{table}[htbp]
\centering
\begin{tabular}{ccc}
\toprule
interval & mesh intervals & runtime (s) \\
\midrule
$[0,\ 2.0\times 10^5]$ & 47 & 13.87\\
$[0,\ 2.0\times 10^4]$ & 38 & 9.406\\
$[0,\ 2.0\times 10^3]$ & 36 & 6.984\\
$[0,\ 2.0\times 10^2]$ & 30 & 5.609\\
$[0,\ 20]$             & 15 & 3.594\\
\bottomrule
\end{tabular}
\caption{Mesh counts and run times for the nonlinear mass-spring optimal-control
problem on increasingly long time intervals.}
\label{tab:1}
\end{table}

Since no closed-form solution is available for this optimal-control system, we
use a self-convergence check for the longest interval reported in the figure,
$T=200$.  The baseline computation starts from the two-interval initial mesh.  We
then take the adaptive terminal mesh, manually double every final subinterval
once, and rerun the solver on this refined mesh.  Comparing the two solutions on
the same target grid gives
\begin{equation}
    \|y_{\rm base}-y_{\rm ref}\|_{\infty}=1.21\times 10^{-12}.
\end{equation}
The maximum differential residual $\|y'-F(y,t)\|_{\infty}$ on the target grid is
$4.80\times10^{-14}$, and the boundary residual is $4.86\times10^{-16}$.  Thus
the adaptive nonlinear solve is self-converged and satisfies the differential
and boundary equations to near machine precision.

With the default merge safety factor, the mesh can retain mild over-refinement.  In
Figure~\ref{fig:1-1}, for instance, the solution is nearly zero on a long
interior interval, yet a few subintervals remain from earlier Newton refinements.
This has little effect on the total cost in the present examples, but more
aggressive coarsening is a possible improvement.

\subsection{Elliptic Jacobi problem}
This system was used in \cite{starrrokhlin} to demonstrate the effectiveness of integral equation methods for systems.  We revisit it here for a different purpose: the problem has multiple solution branches, and Newton iteration can be sensitive to the initial guess and mesh resolution.  It therefore provides a useful test of the continuation and adaptive components of the present solver.  Jacobi elliptic functions arise in, for example, the motion of a pendulum and satisfy the system
\begin{equation}
\left\{
    \begin{aligned}
        & sn'(x)=cn(x)\cdot dn(x),\\
        & cn'(x)=-sn(x)\cdot dn(x),\\
        & dn'(x)=-m\cdot sn(x)\cdot cn(x),\\
        &sn(0)=0,\quad cn(0)=1,\quad dn(L\cdot K)=1.
    \label{eqn: example 2}
    \end{aligned}
\right.
\end{equation}
where $m\in[0,1]$ is fixed to be $1/2$ in our experiments, and $K=K(m)$ is the complete elliptic integral of the first kind,
\[K=K(m)=\int_0^{\frac{\pi}{2}}\frac{d\theta}{\sqrt{1-m\cdot \sin^2\theta}},\]
$L$ is an even number.

\begin{figure}[H]
    \centering
    \begin{minipage}{0.48\textwidth}
        \centering
        \includegraphics[width=\linewidth]{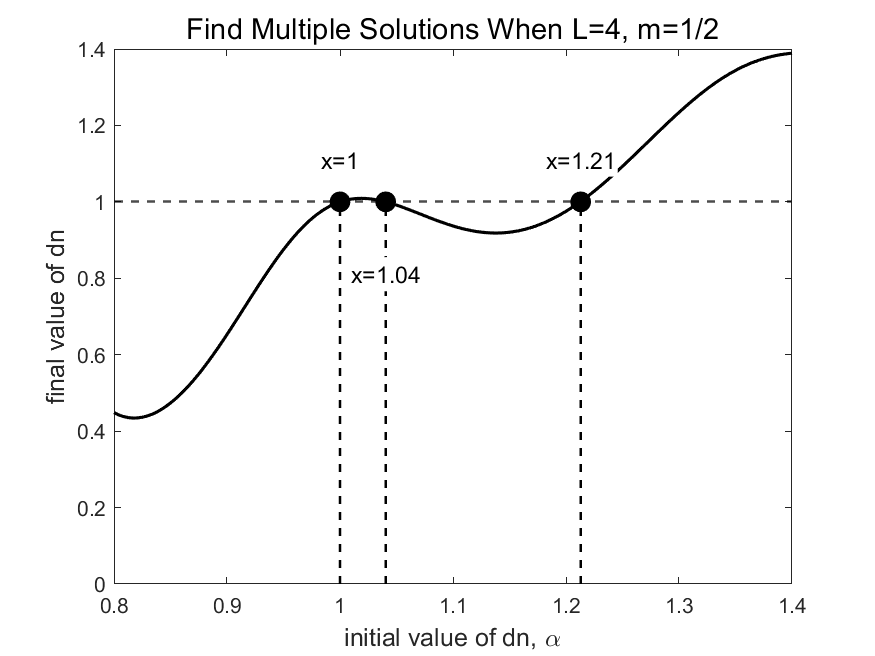}
        \par\smallskip
        {\footnotesize (a) When $L=4$, three solutions are detected.}
    \end{minipage}
    \hfill
    \begin{minipage}{0.48\textwidth}
        \centering
        \includegraphics[width=\linewidth]{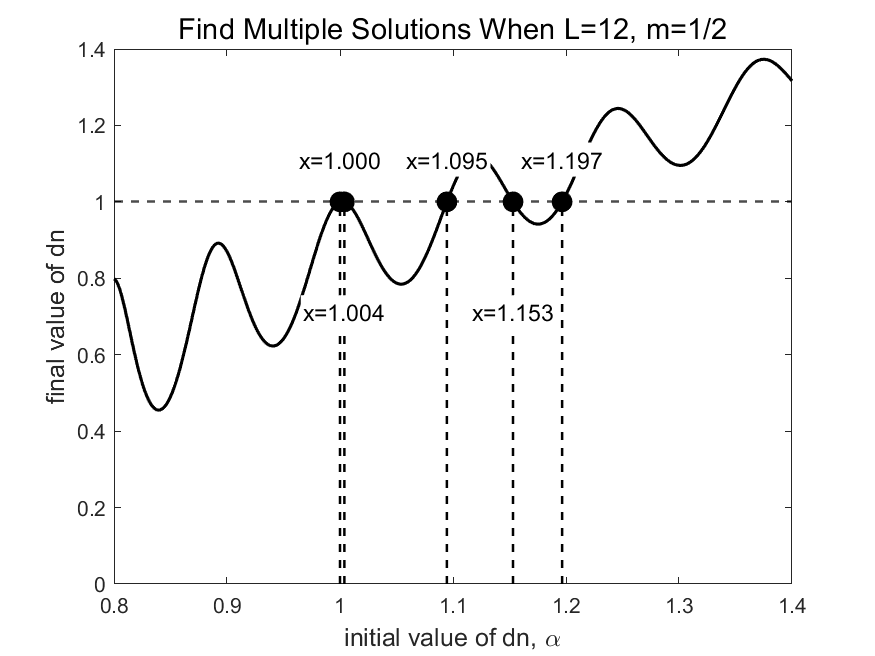}
        \par\smallskip
        {\footnotesize (b) When $L=12$, five solutions are detected.}
    \end{minipage}
    \caption{The plot has $dn(0)$ on the horizontal axis and $dn(1)$ on the vertical axis. A solution to the boundary value problem corresponds to a point where $dn(1) = 1$. Solutions for $L = 4$ and $L = 12$ are displayed in the figure.}
    \label{fig:2-1}
\end{figure}

\begin{figure}[tbp]
    \centering
    \begin{minipage}{0.31\textwidth}
        \centering
        \includegraphics[width=\linewidth]{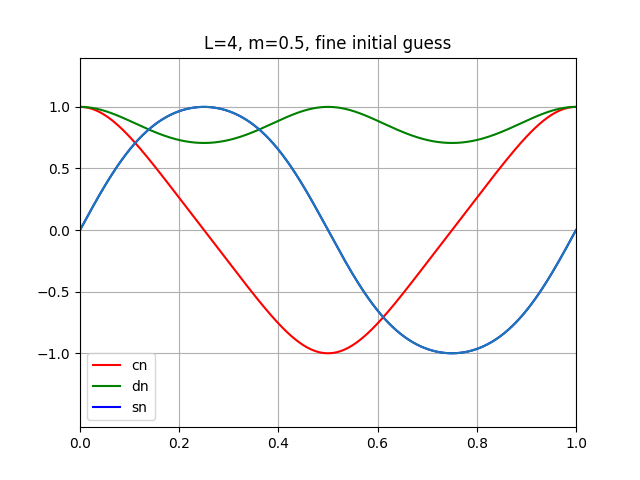}
        \par\smallskip
        {\footnotesize (a) Fine initial guess, $L=4$.}
    \end{minipage}
    \hfill
    \begin{minipage}{0.31\textwidth}
        \centering
        \includegraphics[width=\linewidth]{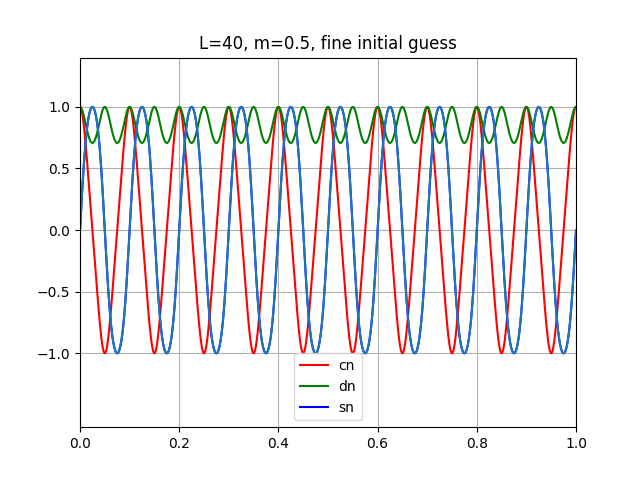}
        \par\smallskip
        {\footnotesize (b) Fine initial guess, $L=40$.}
    \end{minipage}
    \hfill
    \begin{minipage}{0.31\textwidth}
        \centering
        \includegraphics[width=\linewidth]{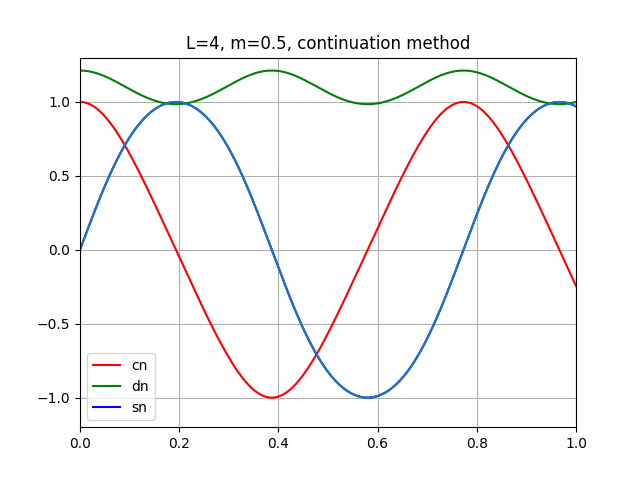}
        \par\smallskip
        {\footnotesize (c) Newton method, $L=4$, $\alpha\approx1.21$.}
    \end{minipage}
    \par\medskip
    \begin{minipage}{0.31\textwidth}
        \centering
        \includegraphics[width=\linewidth]{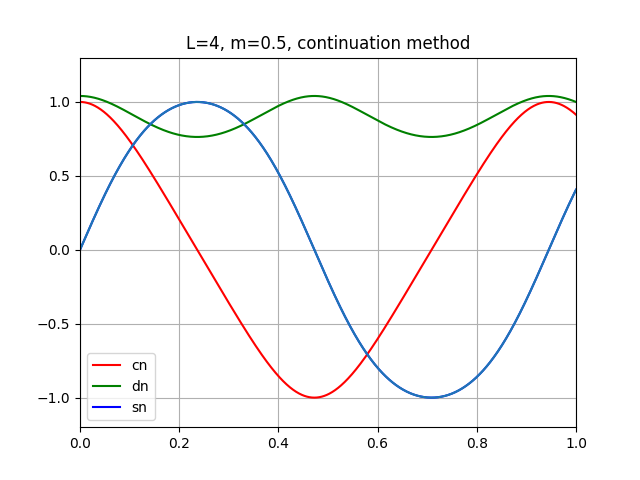}
        \par\smallskip
        {\footnotesize (d) Continuation method, $L=4$, $\alpha\approx1.04$.}
    \end{minipage}
    \hfill
    \begin{minipage}{0.31\textwidth}
        \centering
        \includegraphics[width=\linewidth]{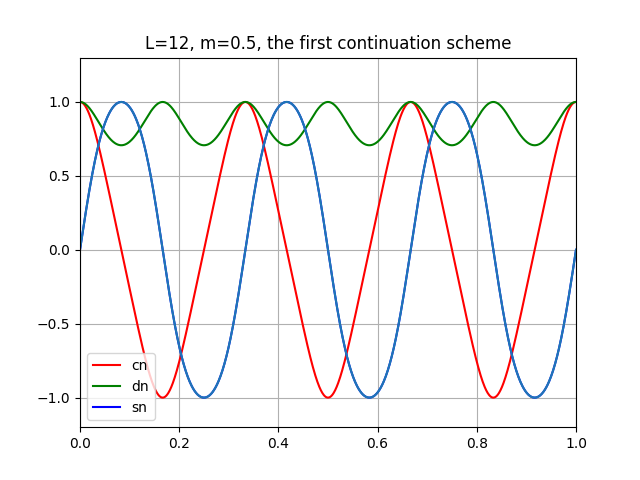}
        \par\smallskip
        {\footnotesize (e) First continuation scheme, $L=12$.}
    \end{minipage}
    \hfill
    \begin{minipage}{0.31\textwidth}
        \centering
        \includegraphics[width=\linewidth]{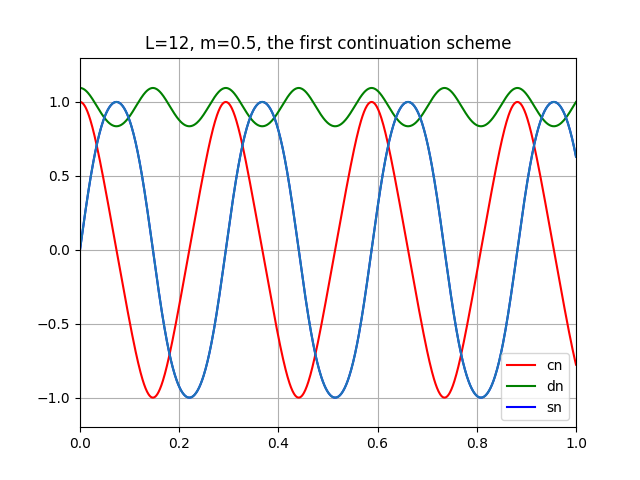}
        \par\smallskip
        {\footnotesize (f) Second continuation scheme, $L=12$.}
    \end{minipage}
    \caption{Numerical solutions for \eqref{eqn: example 2}. All solved in $10^{-12}$ precision. (a): Solution with fine initial guess, $L=4$; (b): Solution with fine initial guess, $L=40$; (c): Solution for $L=4$ without continuation, $\alpha\approx 1.21$; (d): Solution for $L=4$ under continuation method, $\alpha\approx 1.04$; (e): Solution for $L=12$ under the first continuation scheme; (f) Solution for $L=12$ under the second continuation scheme.}
    \label{fig:2-2}
\end{figure}

Following \cite{starrrokhlin}, one natural initial guess is obtained from the case $m=0$, namely $sn(x)=\sin(2\pi x/(4L))$, $cn(x)=\cos(2\pi x/(4L))$, and $dn(x)=0$.  With this initial guess, the present method solves \eqref{eqn: example 2} with $m=1/2$ to about $10^{-12}$ accuracy; see Figure~\ref{fig:2-2}(a),(b).  For larger $L$, Newton's method may fail if the initial mesh is too coarse.  This sensitivity is related to multiple solution branches and a complicated basin of attraction.  To see the source of nonuniqueness, consider the initial value problem

\begin{equation}
\left\{
    \begin{aligned}
        & sn'(x)=cn(x)\cdot dn(x),\\
        & cn'(x)=-sn(x)\cdot dn(x),\\
        & dn'(x)=-m\cdot sn(x)\cdot cn(x),\\
        &sn(0)=0,\quad cn(0)=1,\quad dn(0)=\alpha \quad (\alpha>\sqrt{m})
    \label{eqn: e2IVP}
    \end{aligned}
\right.
\end{equation}
Its solution can be expressed as $\widetilde{sn}(x)=sn_{\widetilde m}(\alpha x)$, $\widetilde{cn}(x)=cn_{\widetilde m}(\alpha x)$, and $\widetilde{dn}(x)=\alpha dn_{\widetilde m}(\alpha x)$, with $\widetilde m=m/\alpha^2$.  Such a solution also satisfies the boundary value problem precisely when the boundary condition on $dn$ is met.  Numerically, additional solutions appear as $L$ grows; see Figure~\ref{fig:2-1}.

We use continuation both to reduce the dependence on the initial guess and to probe different solution branches.  In the experiments below, $p=16$ and $m=1/2$.  For $L=4$, the direct Newton iteration, the continuation method, and the trigonometric initial guess converge to different solutions, as shown in Figure~\ref{fig:2-2}(a),(c),(d).  We do not attempt an exhaustive enumeration of all branches; the point is that different initial guesses and continuation paths can lie in different basins of attraction and therefore lead to distinct solutions, while the adaptive integral-equation solver can still resolve the resulting branches to high accuracy.  This example illustrates why a reliable fast solver for the linearized correction equations is a useful computational foundation for subsequent bifurcation and branch-following calculations.

For $L>4$, the trivial initial guess does not converge without continuation.  To obtain solutions for $L=12$, we tested two continuation schemes:
\begin{itemize}
    \item In the first approach, $m$ is continued from an initial value of 0 to 0.5 while holding $K$ fixed at $K=K(\frac{1}{2})$;
    \item In the second approach, the same continuation on $m$ is performed with $K$ varying, $K=K(m)$.
\end{itemize}
Results are shown in Figure \ref{fig:2-2}, (e) and (f).

\subsection{Troesch equation}
Troesch's problem is a classical nonlinear two-point boundary value problem
with a sharp right boundary layer for large positive $\lambda$:
\begin{equation}
    u''(x)=\lambda\sinh(\lambda u(x)),\qquad
    u(0)=0,
    \quad
    u(1)=1,
    \quad x\in[0,1].
    \label{eqn:troesch_bvp}
\end{equation}
We use the first-order formulation
\begin{equation}
    \phi_1'=\phi_2,
    \qquad
    \phi_2'=\lambda\sinh(\lambda\phi_1),
    \qquad
    \phi_1(0)=0,
    \quad
    \phi_1(1)=1.
    \label{eqn: TroeschBVP}
\end{equation}
The large-parameter regime is a demanding test for nonlinear BVP solvers.  The
regular solution develops an increasingly thin layer near $x=1$, and direct use
of standard collocation solvers such as MATLAB's \texttt{bvp4c} with simple
initial guesses is not reliable in this regime.  Large-parameter Troesch
computations are therefore often treated with carefully designed continuation,
mesh selection, or problem-specific techniques; see, for example,
\cite{Chang2010,Hector2012,Youssef2019}.

For this problem an analytic representation in terms of elliptic functions is
available and can be evaluated to high precision.  Let $s=u_*'(0)>0$, $a=s/2$,
and $m=1-a^2$.  The first integral is
\begin{equation}
    u_*'(x)^2=s^2+4\sinh^2\!\left(\frac{\lambda u_*(x)}{2}\right).
    \label{eqn:troesch_first_integral}
\end{equation}
The regular branch is given by
\begin{equation}
    u_*(x)=\frac{2}{\lambda}\operatorname{arsinh}
       \!\left[a\,\operatorname{sc}(\lambda x\mid m)\right],
    \qquad
    u_*'(x)=\frac{s}{\operatorname{cn}(\lambda x\mid m)},
    \label{eqn:troesch_reference}
\end{equation}
where $\operatorname{sc}=\operatorname{sn}/\operatorname{cn}$ and $m$ is the
elliptic parameter.  The boundary condition $u_*(1)=1$ determines $a$ through
\begin{equation}
    F\!\left(\arctan\!\frac{\sinh(\lambda/2)}{a}
       \,\middle|\,1-a^2\right)=\lambda,
    \label{eqn: ini_tro}
\end{equation}
where $F$ is the incomplete elliptic integral of the first kind.  In the
numerical comparisons below this scalar equation is solved in $\log a$, and the
reference solution is evaluated with \texttt{mpmath} using 100 decimal digits of
working precision.  Recomputing selected values at higher precision confirms at
least 50 correct digits in the reference data.

Let $\widehat u$ denote the numerical solution and $u_*$ the high-precision
reference.  For $f=\widehat u$ or $\widehat u'$ and $f_*=u_*$ or $u_*'$, we
report
\begin{equation}
    E_{2,\mathrm{rel}}(f)
       =\frac{\left(\int_0^1|f(x)-f_*(x)|^2\,dx\right)^{1/2}}
              {\left(\int_0^1|f_*(x)|^2\,dx\right)^{1/2}},
    \qquad
    E_{\infty,\mathrm{mix}}(f)
       =\max_{x_i}\frac{|f(x_i)-f_*(x_i)|}{1+|f_*(x_i)|}.
\label{eqn: tro_metrics}
\end{equation}
The maximum in the mixed norm is sampled on 10,000 equally spaced points; the
term $1+|f_*(x_i)|$ avoids an artificial blow-up where the exact solution is
small.  We also report the absolute derivative error at the left endpoint and
the relative derivative error at the right endpoint,
\begin{equation}
    E_{0,\mathrm{abs}}=|v_h(0)-u_*'(0)|,
    \qquad
    E_{1,\mathrm{rel}}=
       \frac{|v_h(1)-u_*'(1)|}{|u_*'(1)|},
    \label{eqn:troesch_endpoint_errors}
\end{equation}
where $v_h$ is the second component of the computed first-order state.

We first compare direct Newton iteration at $\lambda=20$ with a simple
continuation path $10\to20$, using tolerance $10^{-8}$ throughout.  Both runs
satisfy the internal Newton stopping criterion, but only the continuation run
reaches the high-precision reference to nearly full accuracy; see
Table~\ref{tab:troesch_continuation}.

\begin{table}[tbp]
    \centering
    \caption{Effect of continuation for $\lambda=20$ with tolerance $10^{-8}$. Maximum errors use the 10{,}000-point grid in \eqref{eqn: tro_metrics}; endpoint errors are defined in \eqref{eqn:troesch_endpoint_errors}.}
    \label{tab:troesch_continuation}
    \resizebox{\linewidth}{!}{%
    \begin{tabular}{@{}lrrrrrr@{}}
        \hline
        Method & $E_{2,\mathrm{rel}}(u)$ & $E_{2,\mathrm{rel}}(u')$
        & $E_{\infty,\mathrm{mix}}(u)$ & $E_{\infty,\mathrm{mix}}(u')$
        & $E_{0,\mathrm{abs}}$ & $E_{1,\mathrm{rel}}$ \\
        \hline
        Direct Newton & $9.69743\times10^{-3}$ & $1.22246\times10^{-4}$ & $1.35894\times10^{-3}$ & $1.03871\times10^{-2}$ & $2.94783\times10^{-9}$ & $1.85253\times10^{-10}$ \\
        $10\to20$ & $1.31244\times10^{-10}$ & $2.04517\times10^{-11}$ & $3.79365\times10^{-11}$ & $3.01909\times10^{-10}$ & $2.38731\times10^{-13}$ & $1.24368\times10^{-12}$ \\
        \hline
    \end{tabular}}
\end{table}

The direct run no longer behaves reliably when the tolerance is tightened to
$10^{-9}$, indicating convergence to a nearby numerical profile rather than to
the regular branch at the requested accuracy.  The continuation path supplies a
more accurate initial state for the target parameter and removes this failure.

Table~\ref{tab:troesch_accuracy} reports the four tests on $\lambda = 10,\, 20,\, 30,\,$ and $40$. We use direct Newton at $\lambda=10$ and continuation on others; each sequence uses the target tolerance listed in the table at every stage.

\begin{table}[tbp]
    \centering
    \caption{Reference-solution errors for the four Troesch tests. The derivative approximation is the second physical state component $v_h$. The last two columns give the absolute derivative error at $x=0$ and relative derivative error at $x=1$.}
    \label{tab:troesch_accuracy}
    \resizebox{\linewidth}{!}{%
    \begin{tabular}{@{}rrrrrrrr@{}}
        \hline
        $\lambda$ & Tolerance & $E_{2,\mathrm{rel}}(u)$ & $E_{2,\mathrm{rel}}(u')$
        & $E_{\infty,\mathrm{mix}}(u)$ & $E_{\infty,\mathrm{mix}}(u')$
        & $E_{0,\mathrm{abs}}$ & $E_{1,\mathrm{rel}}$ \\
        \hline
        10 & $10^{-12}$ & $2.92607\times10^{-15}$ & $8.32579\times10^{-16}$ & $6.00109\times10^{-16}$ & $2.13348\times10^{-15}$ & $5.93118\times10^{-16}$ & $1.51479\times10^{-16}$ \\
        20 & $10^{-9}$ & $4.42617\times10^{-12}$ & $2.45051\times10^{-13}$ & $3.91681\times10^{-13}$ & $1.21461\times10^{-12}$ & $2.39450\times10^{-13}$ & $1.30967\times10^{-13}$ \\
        30 & $10^{-5}$ & $1.65600\times10^{-7}$ & $9.95677\times10^{-7}$ & $1.04130\times10^{-7}$ & $1.98420\times10^{-6}$ & $1.00087\times10^{-13}$ & $1.01682\times10^{-8}$ \\
        40 & $10^{-4}$ & $5.91342\times10^{-8}$ & $2.84482\times10^{-6}$ & $3.05401\times10^{-9}$ & $7.97122\times10^{-8}$ & $1.15709\times10^{-14}$ & $2.21180\times10^{-8}$ \\
        \hline
    \end{tabular}}
\end{table}

For $\lambda=30$ and $40$, the final adaptive meshes contain 162 and 1727 subintervals, respectively, with smallest lengths approximately $2^{-32}$ and $2^{-46}$.  The hardest case therefore drives the finest spacing close to the limits of double precision.  The deterioration of the forward error in this regime is consistent with the severe conditioning of the large-parameter Troesch problem rather than with a failure to resolve the boundary layer.

\subsection{Nonlinear Singularly Perturbed BVP}
\begin{figure}
    \centering
    \includegraphics[width=0.8\linewidth]{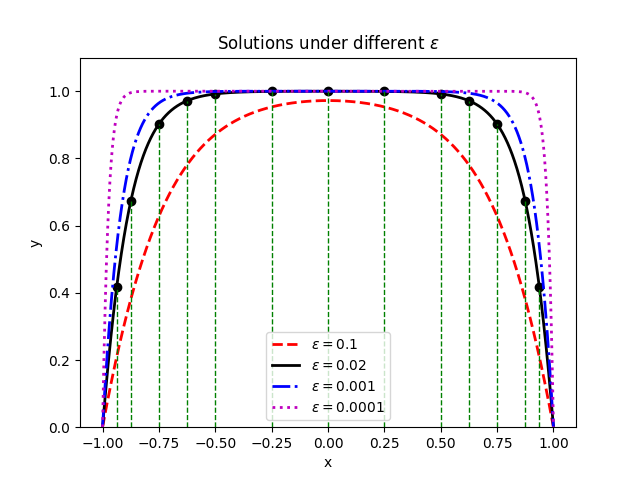}
    \caption{Solutions of the nonlinear singularly perturbed problem \eqref{eqn: PerBVP} for different values of $\epsilon$, computed by continuation.}
    \label{fig:solper}
\end{figure}

The final example is a nonlinear singularly perturbed BVP with boundary layers near both endpoints.  The scalar equation is
\begin{equation}
    \epsilon y''=1-y^2,
    \qquad -1\le x\le1,
    \qquad y(-1)=y(1)=0.
\end{equation}
We write it as the first-order system
\begin{equation}
\left\{
    \begin{aligned}
        & \phi_1'=\phi_2\\
        &\phi_2'= (\phi_1^2-1)/\epsilon
    \end{aligned}
\right.
\label{eqn: PerBVP}
\end{equation}
with boundary condition $\phi_1(-1)=\phi_1(1)=0$.  As $\epsilon$ decreases, the solution develops increasingly sharp endpoint transitions.  Starting Newton's method from the zero initial guess fails for small $\epsilon$, so we use continuation in $\epsilon$.  Figure~\ref{fig:solper} shows the solutions for $\epsilon=0.1,0.02,0.001,0.0001$, computed sequentially to a tolerance of $10^{-12}$.  The vertical lines show the adaptive mesh for one representative continuation stage; Table~\ref{tab:singular_grid_evolution} reports the terminal grid after each continuation solve in the current run.

\begin{table}[htbp]
\centering
\begin{tabular}{rrrr}
\toprule
continuation stage & $\epsilon$ & intervals & $h_{\min}$ \\
\midrule
1 & $10^{-1}$ & 4  & $5.00\times10^{-1}$ \\
2 & $2\times10^{-2}$ & 6  & $2.50\times10^{-1}$ \\
3 & $10^{-3}$ & 10 & $6.25\times10^{-2}$ \\
4 & $10^{-4}$ & 14 & $1.56\times10^{-2}$ \\
\bottomrule
\end{tabular}
\caption{Evolution of the terminal adaptive mesh along the continuation path for
the nonlinear singular perturbation example.  The largest interval remains
$h_{\max}=0.5$ in these runs, while the smallest interval decreases as the
endpoint layers sharpen.}
\label{tab:singular_grid_evolution}
\end{table}

This example also illustrates the interaction between continuation and adaptivity.  Continuation reduces the number of Newton iterations required to reach the final state and therefore avoids the accumulation of unnecessary refinements from poor intermediate iterates.  In this sense, a suitable continuation path improves both nonlinear convergence and mesh selection.

\section{Conclusion}\label{sec:conclusion}
In this paper, we have presented a fast automatically adaptive solver for nonlinear two-point boundary value systems.  The main idea is to formulate the nonlinear problem as a nonlinear second-kind integral equation and to solve each Newton or continuation correction with a fast direct solver for the corresponding linear second-kind equation.  A central message is that progress on fast direct integral-equation solvers for linear BVPs is not confined to linear problems: once the boundary formulation is completed and the linear solver is embedded in Newton, continuation, and adaptive refinement loops, its favorable properties are inherited by a nonlinear solver.  The resulting method retains the well-conditioned nature of integral equations, high-order local Chebyshev discretization, linear complexity in the number of mesh points for fixed system dimension and local order, and automatic concentration of degrees of freedom near localized features.

The method is designed for the difficult regimes that motivate adaptive BVP solvers: boundary layers, rapid oscillations, long intervals with localized dynamics, singular perturbations, hypersensitive optimal-control systems, and problems whose boundary conditions are degenerate for the simplest background Green's function.  The boundary-condition preprocessing makes the formulation applicable to general first-order systems, including high-order scalar equations after conversion to system form.  The numerical examples demonstrate high accuracy with strongly nonuniform meshes selected automatically.

The next step is to use the solver as a computational kernel for bifurcation and solution-branch exploration.  The continuation framework already separates the path-following logic from the fast linearized solves, which makes it suitable for detecting multiple solutions and following branches through parameter space.  Further work will focus on robust Newton globalization, branch switching, and eventually a fast robust method for the computation of folds and solution branches.

\section*{Acknowledgements}
We would like to thank Leslie Greengard, Charles Epstein, and Manas Rachh for many helpful discussions. The authors were supported in part by National Key Research and Development Program of China under grant 2023YFA1008902, and in part by the National Natural Science Foundation of China under grant 12301515.

\bibliographystyle{plain}
\bibliography{bvpnonlin}

\end{document}